\documentclass[a4paper,12pt]{article}        
              
\usepackage{amsthm,amsfonts,amsmath,amssymb} 
\usepackage{mathrsfs}                        
\usepackage{dsfont}    
\usepackage{hyperref} 
\usepackage[T1]{fontenc}     

\usepackage[margin=7em]{geometry}                  

\newcommand{\ud}{\mathrm{d}}

\newcommand{\e}{\mathrm{e}}
\newcommand{\CR}{\mathds{R}}
\newcommand{\CZ}{\mathds{Z}}
\newcommand{\CC}{\mathds{C}}

\newcommand{\dep}[2]{\displaystyle\frac{\partial\emph{$#1$}}{\partial{
\emph{$#2$}}}}
\newcommand{\depp}[2]{\frac{\partial\emph{$#1$}}{\partial{
\emph{$#2$}}}}
\newcommand{\dev}[2]{\displaystyle\frac{\ud\emph{$#1$}}{\ud\emph{$#2$}}}
\newcommand{\devv}[2]{\frac{\ud\emph{$#1$}}{\ud\emph{$#2$}}}
\newcommand{\fraction}{\displaystyle\frac}
\newcommand{\integral}[4]{\displaystyle\int_{\emph{$#1$}}^{\emph{$#2$}}
{\emph{$#3$}} \ \ud\emph{$#4$}}

\newcommand{\somatorio}[2]{\sum\limits_{\emph{$#1$}}^{\emph{$#2$}}}
\newcommand{\limite}[2]{\displaystyle\lim_{{\emph{$#1$}}\to{\emph{$#2$}}}}

\theoremstyle{plain} 
\newtheorem{teo}{Theorem}[section]
\newtheorem{defi}{Definition}[section]
\newtheorem{lem}{Lemma}[section]
\newtheorem{prop}{Proposition}[section]
\newtheorem{coro}{Corollary}[section]

\theoremstyle{definition}
\newtheorem{exam}{Example}[section]
\newtheorem{rema}{Remark}[section]

\begin{document}

\title{Circle actions and geometric quantisation}
\author{Romero Solha\footnote{Departamento de Matemática, Universidade Federal de Minas Gerais. Postal address: Avenida Antonio Carlos, 6627 - Caixa Postal 702 - CEP 30161970 - Belo Horizonte, MG. Email address: romerosolha@gmail.com } \thanks{Partially supported by Erasmus Mundus External Cooperation Window EU-Brazil Startup 2009-2010 project; Contact And Symplectic Topology, Research Networking Programme of the European Science Foundation; Geometry, Mechanics and Control Theory network; DGICYT/FEDER project MTM2009-07594: Estructuras Geometricas: Deformaciones, Singularidades y Geometria Integral; MINECO project MTM2012-38122-C03-01: Geometria Algebraica, Simplectica, Aritmetica y Aplicaciones; CAPES/PNPD.}}

\date{\today}

\maketitle


\begin{abstract}
The aim of this article is to present unifying proofs for results in geometric quantisation with real polarisations by exploring the existence of symplectic circle actions. It provides an extension of Rawnsley's results on the Kostant complex, and gives an alternative proof for \'{S}niatycki's and Hamilton's theorems; as well as, a partial result for the focus-focus contribution to geometric quantisation. 
\end{abstract} 



\section{Introduction}

\hspace{1.5em}Geometric quantisation tries to associate a Hilbert space to a symplectic manifold via a complex line bundle. Although it is possible to describe the canonical quantisation using this language, most of the difficulties arise when one tries to mimic this procedure for symplectic manifolds which are not naturally cotangent bundles. Those appear in the context of reduction and are far from being artificial mathematical models.

The first difficulty is to isolate in a global way position and momentum, in order to define wave functions from sections of a complex line bundle over the symplectic manifold. This is done by introducing polarisations, which, roughly speaking, are lagrangian foliations. The second issue, that will not be addressed here, is how to define a Hilbert structure; however, all examples treated in this article have a natural one.

Usually, the quantum phase space is constructed using global sections of the line bundle which are flat along the polarisation. In case these global sections do not exist, Kostant suggested to associate quantum states to elements of higher cohomology groups, and to built the quantum phase space from these groups: by considering cohomology with coefficients in the sheaf of flat sections.

At least two approaches can be used to compute these cohomology groups: \v{C}ech and de Rham. The results of Hamilton \cite{Ha} and Hamilton and Miranda \cite{HaMi} are based on a \v{C}ech approach, this article takes the de Rham point of view, by finding a resolution for the sheaf. Following Kostant \cite{Sni,JHR}, a resolution for the sheaf of flat sections can be obtained by twisting the sheaves relative to the foliated complex induced by the polarisation with the sheaf of flat sections.

This article follows closely Rawnsley's ideas \cite{JHR} and explores the existence of circle actions, in the particular case of real polarisations, to provide a different proof for the theorems of \'{S}niatycki \cite{Sni} and Hamilton \cite{Ha}. The tools developed here highlight and unravel the role played by symplectic circle actions in known results in geometric quantisation. Not only that, this approach casts some light on a conjecture about the contributions coming from focus-focus type of singularities. 

The rest of the article is organised as follows. Section 2 introduces the basic definitions of geometric quantisation. Section 3 summarises relevant results about a resolution of the sheaf of flat sections. Section 4 explores the existence of circle actions: it further develops results from \cite{JHR} and it contains the main tools of this article. Section 5 presents the important notion of Bohr-Sommerfeld leaves. In section 6 the tools developed in section 4 are used to compute the geometric quantisation of local and semilocal models, in particular near a focus-focus singularity and fibre. Finally, using the results of section 4 and 6, section 7 provides an alternative proof for \'{S}niatycki's \cite{Sni} and Hamilton's \cite{Ha} theorem.\par

Throughout this article and otherwise stated, all the objects considered will be $C^\infty$; manifolds are real, Hausdorff, paracompact, and connected; $C^\infty(V)$ denotes the set of complex-valued functions over $V$; and the units are such that $\hbar=1$. 


\subsection{Acknowledgements}

\hspace{1.5em}There is no way to overemphasise the importance of Eva Miranda in this work: she introduced the author of this article to the subject, and also read and commented on drafts of the article. The author also wants to express his gratitude to Francisco Presas and an anonymous referee for their valuable observations and comments. 


\section{Geometric quantisation}


\subsection{Prequantisation}\label{prequantum}

\hspace{1.5em}This subsection deals with some concepts needed to define quantum states. The first attempt was to see them as sections of a complex line bundle over the symplectic manifold, the so-called prequantum line bundle. The other notion described here, polarisation, is a way to define a global distinction between momentum and position.

Using an isomorphism between the \v{C}ech cohomology $\check{H}^2(M;\CR)$ and de Rham cohomology $H^2_{dR}(M;\CR)$, a closed $2$-form is integral if and only if it is in the image of the homomorphism between $\check{H}^2(M;\CZ)$ and $\check{H}^2(M;\CR)$:  
\begin{center}
\begin{tabular}{c c c c c}
$\CR$ &  & $\check{H}^2(M;\CR)$ & $\longleftrightarrow$ & $H^2_{dR}(M;\CR)$ \\
$\uparrow$ & & $\uparrow$ & & \\
$\CZ$ & & $\check{H}^2(M;\CZ)$  & & 
\end{tabular} 
\end{center}
the homomorphism between $\check{H}^2(M;\CZ)$ and $\check{H}^2(M;\CR)$ is induced by a homomorphism between $\CZ$ and $\CR$.

\begin{defi} A symplectic manifold $(M,\omega)$ such that the de Rham class $[\omega]$ is integral is called prequantisable.
\end{defi}

\begin{defi} A prequantum line bundle of $(M,\omega)$ is a hermitian line bundle over $M$ with connexion, compatible with the hermitian structure, $(L,\nabla^{\omega})$ that satisfies $curv(\nabla^{\omega})=-i\omega$.
\end{defi}

\begin{exam} Any exact symplectic manifold satisfies $[\omega]=0$, in particular cotangent bundles with the canonical symplectic structure. The trivial line bundle is an example of a prequantum line bundle in this case.\hfill $\Diamond$\par\vspace{0.5em}
\end{exam}

The following theorem (a proof can be found in \cite{Ko}) provides a relation between the above definitions:

\begin{teo} A symplectic manifold $(M,\omega)$ admits a prequantum line bundle $(L,\nabla^{\omega})$ if and only if it is prequantisable.
\end{teo}

When the symplectic manifold is exact its symplectic form belongs to the kernel of the homomorphism $\check{H}^2(M;\CR)\to H^2_{dR}(M;\CR)$, and the associated prequantum line bundles  are flat hermitian line bundles. Up to bundle isomorphisms, a flat hermitian line bundle is determined by a homomorphism between the fundamental group of $M$, $\pi_1(M)$, and the unitary group, $U(1)$: the holonomy of a flat connexion.   

\begin{lem}\label{flatbundle}At a submanifold $N\subset M$ where $curv(\nabla^\omega)\big{|}_{TN}=-i\ud\Theta$, $(L,\nabla^\omega)\big{|}_{N}\cong (E,\nabla^0)\otimes (\CC\times N,\nabla^\Theta)$ with $(E,\nabla^0)$ being a flat hermitian line bundle over $N$ and $curv(\nabla^\Theta)=-i\ud\Theta$.  
\end{lem}

Classically, a real polarisation $\mathcal{F}$ is an integrable subbundle of $TM$ (the bundle $T\mathcal{F}$) whose leaves are lagrangian submanifolds: i.e. $\mathcal{F}$ is a lagrangian foliation. But due to the example below another definition is considered.

An integrable system on a symplectic manifold $(M,\omega)$ of dimension $2n$ is a set of $n$ real-valued functions, $f_1,\dots,f_n\in C^\infty(M)$, satisfying $\ud f_1\wedge\cdots\wedge\ud f_n\neq 0$ over an open dense subset of $M$ and $\{f_j,f_k\}_\omega=0$ for all $j,k$. The mapping $F=(f_1,\dots,f_n):M\to\CR^n$ is called a moment map.

The Poisson bracket is defined by $\{f,g\}_\omega=X_f(g)$, where $X_f$ is the unique vector field defined by the equation $\imath_{X_f}\omega=-\ud f$, called the hamiltonian vector field of $f$.

The distribution generated by the hamiltonian vector fields of the moment map is involutive because $[X_f, X_g]=X_{\{f,g\}_\omega}$. Since $0=\{f_j,f_k\}_\omega=\omega(X_{f_j}, X_{f_k})$, the leaves of the associated (possibly singular) foliation are isotropic submanifolds and they are lagrangian at points where the functions are functionally independent. This is an example of a generalised real polarisation ---i.e. an integrable distribution on $TM$ whose leaves are lagrangian submanifolds, except for some singular leaves.

\begin{defi}\label{polarisationdef} A real polarisation $\mathcal{P}$ is an integrable (in the Sussmann's sense \cite{SSNN}) distribution of $TM$ whose leaves are generically lagrangian. The complexification of $\mathcal{P}$ is denoted by $P$ and will be called polarisation.
\end{defi}

The most relevant polarisation for this article is $\langle X_{f_1},...,X_{f_n}\rangle_{C^\infty(M)}$: the distribution of the hamiltonian vector fields $X_{f_i}$ of the components $f_i$ of an integrable system $F=(f_1,\dots,f_n):M\to\CR^{n}$.

\begin{defi}
A point $p\in M$ in a $2n$ dimensional symplectic manifold $(M,\omega)$ is a nondegenerate singular point of  Williamson type $(k_e,k_h,k_f)$ of an integrable system $F=(f_1,\dots,f_n):M\to\CR^{n}$ if $p$ is a critical point of rank $n-k_e-k_h-2k_f$ and, in some symplectic system of coordinates $(x_1,y_1,\dots,x_n,y_n)$ where $\omega=\somatorio{i=1}{n} \ud x_i\wedge \ud y_i$, the quadratic parts of $f_1,\dots,f_n$ can be written as:
  \begin{equation}
    \begin{array}{lcr}
      h_i = x_i^2 + y_i^2 & \textup{for }  1 \leq i \leq k_e \ , &
\textup{(elliptic)}  \\
      h_i = x_iy_i  &  \textup{for }  k_e+1 \leq i \leq k_e+k_h \ , &
      \textup{(hyperbolic)}\\
      \begin{cases}
        h_i = x_i y_i + x_{i+1} y_{i+1} , \\
        h_{i+1} = x_i y_{i+1}-x_{i+1} y_i
      \end{cases} &
      \begin{array}{c}
        \textup{for }  i = k_e+k_h+ 2j-1,\\ 1\leq j \leq k_f
      \end{array}
      &
      \textup{(focus-focus pair)}
 \end{array}
  \end{equation}
 
\end{defi}

\begin{exam} For the simple pendulum the stable equilibrium point is an elliptic singularity, whilst the unstable one is a hyperbolic. The spherical pendulum has a stable equilibrium point that is a purely elliptic singularity. The unstable equilibrium point is a focus-focus singularity.\hfill $\Diamond$\par\vspace{0.5em}
\end{exam}

Here is an example of a real polarisation that do not come from an integrable system.

\begin{exam} The action of $S^1$ on $S^1\times S^1$ given by $(z,x,y)\mapsto (z\cdot x,y)$, with $z,x,y\in S^1$, is symplectic (taking as symplectic form the area form of the torus). Because there are no fixed points, this action cannot be hamiltonian ---otherwise, one would have a function over a compact manifold without critical points.\hfill $\Diamond$\par\vspace{0.5em}
\end{exam}

Henceforth, $(L,\nabla^{\omega})$ will be a prequantum line bundle and $P$ the complexification of a real polarisation $\mathcal{P}$ of $(M,\omega)$.


\subsection{Geometric quantisation à la Kostant}\label{gqkostant}

\hspace{1.5em}The original idea of geometric quantisation is to associate a Hilbert space to a symplectic manifold via a prequantum line bundle and a polarisation. Usually this is done using global flat sections of the line bundle; in case these global sections do not exist, one can define geometric quantisation via higher cohomology groups by considering cohomology with coefficients in the sheaf of flat sections.

The existence of global nonzero flat sections is a nontrivial matter, even when $M$ is not compact. Actually, Rawnsley \cite{JHR} (also proposition \ref{rawaction} in this article, under slightly different hypotheses) showed that the existence of a $S^1$-action may be an obstruction for the existence of nonzero global flat sections.

The cotangent bundle of the circle, endowed with the canonical symplectic structure and a polarisation induced by a circle action, provides an explicit example of the nonexistence of nonzero global flat sections.

\begin{exam}\label{cilindro}
Consider $M=\CR\times S^1$ with coordinates $(x,y)$ and $\omega=\ud x\wedge\ud y$. Take as $L$ the trivial complex line bundle with connexion 1-form $\Theta=x\ud y$, with respect to the unitary section\footnote{Sections of $L$ can be represented by complex-valued functions over local trivialisations. When there is an identification between a section $s$ of the line bundle and a complex-valued function $f$, the bundle isomorphism will be omitted, for the sake of simplicity, and the equality $s=f$ will be used.} $\e^{ix}$, and  $P=\left\langle\depp{}{y}\right\rangle_{C^\infty(M)}$. 

Flat sections, $s(x,y)=f(x,y)\e^{ix}$, satisfy 
\begin{equation}
\left[\nabla^\omega_{\depp{}{y}}s\right](x,y)=\left(\dep{f}{y}(x,y)-ixf(x,y)\right)\e^{ix}=0 \ .
\end{equation}Thus, $s(x,y)=g(x)\e^{ixy}\e^{ix}$, for some function $g$, and it has period $2\pi$ in $y$ if and only if $x\in\CZ$, for $S^1$ the unity circle: flat sections are only well-defined for the set of points with $x\in\CZ$. \hfill $\Diamond$\par\vspace{0.5em}
\end{exam}

Global flat sections are, then, absent generally, and if one insists on using flat sections as analogue for quantum states, they is forced to work with distributions with support over Bohr-Sommerfeld leaves (definition \ref{B-Sleaf}), or deal with sheaves and higher order cohomology groups. In this article, only the sheaf approach is treated: as suggested by Kostant. 

\begin{defi} Let $\mathcal{J}$ denotes the sheaf of sections of a prequantum line bundle $L$ such that for each open set $V\subset M$ the set $\mathcal{J}(V)$ is the module, over the ring of complex-valued functions of $V$ which are constant along the leaves of a polarisation $P$, of sections of $L$ defined over $V$ satisfying $\nabla^{\omega}_X s=0$ for all vector fields $X$, defined over $V$, of $P$. $\mathcal{J}$ is called the sheaf of flat sections.  
\end{defi}
 
Consider the triplet: prequantisable symplectic manifold $(M,\omega)$, prequantum line bundle $(L,\nabla^{\omega})$, and polarisation $P$.
 
\begin{defi} The quantisation of $(M,\omega,L,\nabla^{\omega},P)$ is given by
\begin{equation}
\mathcal{Q}(M)=\displaystyle\bigoplus_{k\geq 0}\check{H}^k(M;\mathcal{J}) \ ,
\end{equation}where $\check{H}^k(M;\mathcal{J})$ are \v{C}ech cohomology groups with values in the sheaf $\mathcal{J}$. In this case, one implicitly assumes the extra structures and calls $M$ a quantisable manifold.  
\end{defi} 

The present article can be summarised as an approach to compute and understand the features of these cohomology groups when there is a symplectic circle action preserving the polarisation.

\begin{rema} Even though $\mathcal{Q}(M)$ is just a vector space and a priori has no Hilbert structure, it will be called quantisation. The true quantisation shall be the completion of the vector space $\mathcal{Q}(M)$, after a Hilbert structure is given, together with a Lie algebra homomorphism (possibly defined over a smaller subset) between the Poisson algebra of $C^\infty(M)$ and operators on the Hilbert space. In spite of the problems that may exist in order to define geometric quantisation using $\mathcal{Q}(M)$, the first step is to compute this vector space.
\end{rema}

\begin{rema} Flat sections behave in a different fashion for Kähler polarisations (e.g. for compact manifolds $\mathcal{Q}(M)$ is finite dimensional). This article does not deal with this case; however, much can be found in the literature (e.g. \cite{GuSt,Ha} and references therein). There is another aspect of the theory that will be left aside by this article: metaplectic correction. To imbue $\check{H}^0(M;\mathcal{J})\cong\{s\in\Gamma(L) \ ; \ \nabla^\omega_X s=0 \ \forall \ X\in P\}$ with a Hilbert structure, Kostant and Blattner \cite{Ko2,Blt} introduced half-forms on geometric quantisation\footnote{It is not clear who did what, but both Kostant and Blattner say that it has roots on a joint work of them with Sternberg.}. Besides inducing an inner product, half-forms also make a correction to the spectrum of the operators (Blattner, Rawnsley, Simms and \'{S}niatycki are referred to for this in \cite{Ko1,JHR,Sni}), this correction does not always behaves as one would like, though (e.g. \cite{DET}). 
\end{rema}


\section{Resolution approach}

\hspace{1.5em}Following Rawnsley \cite{JHR}, given a prequantisable symplectic manifold $(M,\omega)$ with polarisation $P$ and prequantum line bundle $(L,\nabla^{\omega})$, it is possible to construct a fine resolution for the sheaf of flat sections $\mathcal{J}$. Using the results of section \ref{circleraw}, it is even possible to do it when $P$ has nondegenerate singularities: this is the content of theorems \ref{poincaeli3} and \ref{poincafuck}. 

The propositions contained in this section are extensions of the results in \cite{JHR}; it is mainly an opportunity to fix notation, the replacement of a subbundle of $TM$ by an integrable distribution offers no obstruction and, therefore, proofs are omitted.\par 

The set $C^{\infty}(M)$ is a commutative $\CC$-algebra and the polarisation induced by a integrable system is both a $C^{\infty}(M)$-module and a $\CC$-Lie algebra; indeed, a Lie subalgebra of vector fields of $M$. The Lie algebra and $C^\infty(M)$-module structures are compatible in such a way that $(P,C^{\infty}(M),\CC)$ is an example of a Lie pseudoalgebra (see \cite{Mkz} for precise definitions and a nice account for the history and, various, names of this structure).

\begin{defi} Let $\Omega_P^k(M)$ denote  the space of multilinear maps 
\begin{equation*}
\mathrm{Hom}_{C^\infty(M)}(\wedge^k_{C^\infty(M)}P; C^\infty(M)) \ ,
\end{equation*} and $S_P^k(L):=\Omega_P^k(M)\otimes_{C^\infty(M)}\Gamma(L)$. Then, ${\Omega_P}^\bullet(M):=\displaystyle\bigoplus_{k\geq 0}\Omega_P^k(M)$ is the space of polarised forms, and the space of line bundle valued polarised forms is ${S_P}^\bullet(L):=\displaystyle\bigoplus_{k\geq 0}S_P^k(L)$.
\end{defi}

The restriction of the connexion $\nabla^{\omega}$ to the polarisation, 
\begin{equation}
\nabla:=\nabla^{\omega}\big{|}_P:\Gamma(L)\to\Omega_P^1(M)\otimes_{C^\infty(M)}\Gamma (L)
\end{equation}satisfies (by definition) the following property:
\begin{equation}
\nabla(fs)=\ud_Pf\otimes s + f\nabla s\ ,
\end{equation}for any $f\in C^\infty(M)$ and $s\in\Gamma(L)$.

Therefore, $\nabla:S^0_P(L)\to S^1_P(L)$ and ${S_P}^\bullet(L)$ has a module structure which enables an extension of $\nabla$ to a derivation of degree $+1$ on the space of line bundle valued polarised forms, as follows: if $\alpha\in\Omega_P^k(M)$ and $\boldsymbol\beta=\beta\otimes s\in S_P^l(L)$, 
\begin{equation}
\alpha\wedge\boldsymbol\beta=\alpha\wedge(\beta\otimes s):=(\alpha\wedge\beta)\otimes s \ 
\end{equation} defines a left multiplication of the ring ${\Omega_P}^\bullet(M)$ on ${S_P}^\bullet(L)$.

\begin{defi} The derivation on ${S_P}^\bullet(L)$ is given by the degree $+1$ map $\ud^\nabla:{S_P}^\bullet(L)\to {S_P}^\bullet(L)$, 
\begin{equation}
\ud^\nabla(\alpha\otimes s):=\ud_P\alpha\otimes s+(-1)^k\alpha\wedge\nabla s \ .
\end{equation}
\end{defi}

The exterior derivative $\ud_P$ is  the restriction of the de Rham differential to the directions of the polarisation. Namely, given $\alpha\in \Omega_P^k(M)$ and $Y_1,\dots,Y_{k+1}\in P$, it is defined by:
\begin{eqnarray}
\ud_P\alpha(Y_1,\dots,Y_{k+1})&=&\somatorio{i=1}{k+1}(-1)^{i+1}Y_i(\alpha(Y_1,\dots,\hat{Y_i},\dots,Y_{k+1}))  \\
&&+\somatorio{i<j}{}(-1)^{i+j}\alpha([Y_i,Y_j],Y_1,\dots,\hat{Y_i},\dots,\hat{Y_j},\dots,Y_{k+1}) \ .  \nonumber
\end{eqnarray}

\begin{prop}\label{propcurv} If $\alpha\in\Omega_P^k(M)$ and $\boldsymbol\beta\in S_P^l(L)$, then
\begin{equation}
\ud^\nabla(\alpha\wedge\boldsymbol\beta)=\ud_P\alpha\wedge\boldsymbol\beta+(-1)^k\alpha\wedge\ud^\nabla\boldsymbol\beta \ ,
\end{equation}and
\begin{equation}
\ud^\nabla\circ\ud^\nabla\boldsymbol\beta=curv(\nabla^\omega)\big{|}_P\wedge\boldsymbol\beta \ .
\end{equation}
\end{prop}

Since $\omega=i \ curv(\nabla^\omega)$ vanishes along $P$, one has $\ud^\nabla\circ\ud^\nabla=0$.

\begin{coro}\label{cobocoro} $\ud^\nabla$ is a coboundary operator. 
\end{coro}

Corollary \ref{cobocoro} implies that the restriction of the connexion $\nabla^{\omega}$ to the polarisation defines a representation of the Lie pseudoalgebra $(P,C^{\infty}(M),\CC)$ on $\Gamma(L)$: i.e. a Lie algebra representation of $(P, [\cdot,\cdot]\big{|}_P)$ on $\Gamma(L)$ compatible with their $C^{\infty}(M)$-module structures.

\begin{rema}\label{metarema} The only property of $L$ being used here is the existence of flat connexions along $P$; any complex line bundle admitting such a connexion would do, not only a prequantum one ---the results here work if metaplectic correction is included.
\end{rema}

Thus, the associated Lie pseudoalgebra cohomology of this representation, \linebreak $H^\bullet({S_P}^\bullet(L))$, induces a complex (at the sheaf level). If $\mathcal{S}_P^k(L)$ denotes the associated sheaf of $S_P^k(L)$, one can extend $\ud^\nabla$ to a homomorphism of sheaves, $\ud^\nabla:\mathcal{S}_P^k(L)\to\mathcal{S}_P^{k+1}(L)$. $\mathcal{S}_P^0(L)\cong\mathcal{S}$, the sheaf of sections of the line bundle $L$, and $\mathcal{J}$ is isomorphic to the kernel of $\nabla:\mathcal{S}\to\mathcal{S}_P^1(L)$. For $k\neq 0$ and $V\subset M$ any open set, $\mathcal{S}_P^k(L)(V):=\mathrm{Hom}_{C^\infty(V)}(\wedge^k_{C^\infty(V)}P\big{|}_V; \Gamma(L\big{|}_V))$.

\begin{defi} The Kostant complex is 
\begin{equation}\label{eq72} 
0\longrightarrow\mathcal{J}\hookrightarrow\mathcal{S}\stackrel{\nabla}{\longrightarrow}\mathcal{S}_P^1(L)\stackrel{\ud^\nabla}{\longrightarrow}\cdots\stackrel{\ud^\nabla}{\longrightarrow}\mathcal{S}_P^n(L)\stackrel{\ud^\nabla}{\longrightarrow}0 \ .
\end{equation}
\end{defi}

The sheaves $\mathcal{S}^k_P(L)$ are fine: $\Gamma(L)$ and $\Omega^k_P(M)$ are free modules over the ring of functions of $M$, and because of that, they admit partition of unity. Hence, if one can prove a Poincaré lemma, the abstract de Rham theorem implies that the Kostant complex is a fine resolution for $\mathcal{J}$.

There are particular situations in which a Poincaré lemma is available. This is true\footnote{Both \'{S}niatycki and Rawnsley attribute this to Kostant, a proof is provided in \cite{JHR}.} when $\mathcal{P}$ is a subbundle of $TM$, and it can be extended to a more general setting; this article provides Poincaré lemmata when $\mathcal{P}$ has nondegenerate singularities. 

\begin{teo}\label{fineresolution1} The Kostant complex is a fine resolution for $\mathcal{J}$ when $\mathcal{P}$ is a subbundle of $TM$, or when $\mathcal{P}$ is induced by an integrable system whose moment map has only singularities of Williamson type $(k_e,k_h,k_f)$, with $k_e+k_h+2k_f\leq \frac{1}{2}\mathrm{dim}(M)$, and either $k_e\geq 1$ or $k_f\geq 1$. Therefore, each of its cohomology groups, $H^k({S_P}^\bullet(L))$, is isomorphic to $\check{H}^k(M;\mathcal{J})$.
\end{teo}
\textit{Proof:} Theorems \ref{poincaregu}, \ref{poincaeli3}, and \ref{poincafuck} guarantee the existence of Poincaré lemmata, and the abstract de Rham theorem finishes the proof.\hfill $\blacksquare$\par\vspace{0.5em}

Analogously to the de Rham cohomology case, there exists a Mayer-Vietoris sequence for $\mathcal{J}$. One can construct, for each pair of open subsets $V$ and $W$ of $M$, the injective homomorphism 
\begin{equation}
R_{V,W}:\mathcal{S}_P^k(L)(V\cup W)\to\mathcal{S}_P^k(L)(V)\oplus\mathcal{S}_P^k(L)(W)
\end{equation}defined by $R_{V,W}(\boldsymbol{\zeta})=\boldsymbol{\zeta}\big{|}_{TV}\oplus\boldsymbol{\zeta}\big{|}_{TW}$ and the surjective homomorphism 
\begin{equation}
R_{V,V\cap W}-R_{W,V\cap W}:\mathcal{S}_P^k(L)(V)\oplus\mathcal{S}_P^k(L)(W)\to\mathcal{S}_P^k(L)(V\cap W)
\end{equation}defined by $R_{V,V\cap W}-R_{W,V\cap W}(\boldsymbol{\alpha}\oplus\boldsymbol{\beta})=\boldsymbol{\alpha}\big{|}_{T(V\cap W)}-\boldsymbol{\beta}\big{|}_{T(V\cap W)}$. The injectivity of $R_{V,W}$ is due to the local identity property of the sheaves, whilst the surjectivity of  $R_{V,V\cap W}-R_{W,V\cap W}$ comes from the existence of partitions of unity for $\mathcal{S}^k_P(L)$. 

Thanks to the glueing condition of the sheaves, the image of $R_{V,W}$ is equal to the kernel of  $R_{V,V\cap W}-R_{W,V\cap W}$, and the long exact sequence associated to the short exact sequences
\begin{equation} 
0\rightarrow\mathcal{S}_P^k(L)(V\cup W)\hookrightarrow\mathcal{S}_P^k(L)(V)\oplus\mathcal{S}_P^k(L)(W)\twoheadrightarrow\mathcal{S}_P^k(L)(V\cap W)\rightarrow 0 \ 
\end{equation}yields the Mayer-Vietoris sequence for the Kostant complex (or for $\mathcal{J}$, applying theorem \ref{fineresolution1}).

\begin{rema}It will be said that the Mayer-Vietoris argument works if the cohomology groups $H^k({S_P}^\bullet(L))$ vanish when restricted to an intersection $V\cap W$ of open subsets of $M$. This means that $H^k({S_P}^\bullet(L\big{|}_{V\cup W}))\cong H^k({S_P}^\bullet(L\big{|}_V))\oplus H^k({S_P}^\bullet(L\big{|}_W))$. 
\end{rema}   

As expected, the notions of interior product and Lie derivative are available for ${S_P}^\bullet(L)$. The Lie derivative can be seen as a derivation along a flow, but for that, a nontrivial notion of pullback is needed. 

\begin{defi} The contraction between line bundle valued polarised forms and elements of $P$ is given by a map $\boldsymbol{i}:P\times {S_P}^\bullet(L)\to {S_P}^\bullet(L)$ that is a degree -1 map on ${S_P}^\bullet(L)$: i.e. 
\begin{equation}
\boldsymbol{i}_X(\nabla s):=\nabla_X s  \  
\end{equation}and
\begin{equation}
\boldsymbol{i}_X\boldsymbol\beta=\boldsymbol{i}_X(\beta\otimes s):=(\imath_X\beta)\otimes s \ 
\end{equation}hold for each $X\in P$ and $\boldsymbol\beta=\beta\otimes s\in S_P^l(L)$.
\end{defi}

\begin{prop}\label{degreeminus} If $X\in P$, $\alpha\in\Omega_P^k(M)$ and $\boldsymbol\beta\in S_P^l(L)$, then $\boldsymbol{i}_X\circ\boldsymbol{i}_X=0$ and
\begin{equation}
\boldsymbol{i}_X(\alpha\wedge\boldsymbol\beta)=(\imath_X\alpha)\wedge\boldsymbol\beta+(-1)^k\alpha\wedge\boldsymbol{i}_X\boldsymbol\beta \ .
\end{equation}
\end{prop}
 
\begin{defi} The pullback $\boldsymbol{\phi_t}^*$ of $\alpha\otimes s\in S_P^k(L)$ is defined by
\begin{equation}
\boldsymbol{\phi_t}^*(\alpha\otimes s):=(\phi_t^*\alpha)\otimes \Pi_{\phi_t}^{-1}(s\circ\phi_t) \ ;
\end{equation}where, by the bundle automorphism property of the parallel transport,  $\Pi_{\phi_t}^{-1}(s\circ\phi_t)$ denotes the parallel transport between $\phi_t(p)$ and $p$ of $s$ through the integral curve of the flow. 
\end{defi}

\begin{prop}\label{pullcommutes} Let $X\in P$ with flow $\phi_t$, $\alpha\in\Omega_P^k(M)$ and $\boldsymbol\beta\in S_P^l(L)$; then, 
\begin{equation}
\boldsymbol{\phi_t}^*(\alpha\wedge\boldsymbol\beta)=(\phi_t^*\alpha)\wedge\boldsymbol{\phi_t}^*(\boldsymbol\beta) \ ,
\end{equation}and the pullback $\boldsymbol{\phi_t}^*$ commutes with $\ud^\nabla$.
\end{prop}

\begin{defi} The Lie derivative $\pounds^\nabla:P\times {S_P}^\bullet(L)\to {S_P}^\bullet(L)$ is defined by:
\begin{equation}
\pounds^\nabla_X(\boldsymbol\alpha):=\dev{}{t}\boldsymbol{\phi_t}^*\boldsymbol\alpha\Big{|}_{t=0} \ .
\end{equation}
\end{defi}

Cartan's magic formula holds for the Lie derivative on ${S_P}^\bullet(L)$, and it also commutes with the pullback and exterior derivative.

\begin{prop} The Lie derivative $\pounds^\nabla$ commutes with the pullback $\boldsymbol{\phi_t}^*$ and with the derivation $\ud^\nabla$, and it can be characterised by 
\begin{equation}
\pounds^\nabla_X(\boldsymbol\alpha)=\boldsymbol{i}_X\circ\ud^\nabla\boldsymbol\alpha+\ud^\nabla\circ\boldsymbol{i}_X\boldsymbol\alpha \ .
\end{equation}
\end{prop} 
  

\section{Circle actions and homotopy operators}\label{circleraw}

\hspace{1.5em}This section explains the construction of an almost homotopy operator for the Kostant complex when one has a symplectic $S^1$-action, and how this implies the vanishing of the stalks of points with nontrivial holonomy. Most results of this section were previously provided in \cite{JHR} with slightly less general hypothesis; some proofs automatically hold (propositions \ref{almosthomo}, \ref{constholonomy} and \ref{rawaction}), but one (lemma \ref{rawmagic}) had to be adapted.

Let $X\in P$ be a generator of a symplectic $S^1$-action. If $\phi_t$ stands for the flow of $X$ at time $t$, it is possible to define an induced action on $S_P^k(L)$ via $\boldsymbol{\phi_t}^*$. The holonomy of the loop generated by flowing during a time $2\pi$ a point $p\in M$ is an element of $U(1)$ and will be denoted by $hol_{\nabla^\omega}(\gamma)(p)$, and since $\phi_{t+2\pi}=\phi_t$ for every $t\in\CR$:  
\begin{eqnarray*}
\boldsymbol{\phi_{2\pi}}^*(\boldsymbol\alpha)&=&\boldsymbol{\phi_{2\pi}}^*(\alpha\otimes s)=\phi_{2\pi}^*\alpha\otimes\Pi^{-1}_{\phi_{2\pi}}(s\circ\phi_{2\pi}) \\
&=& \alpha\otimes ({{hol_{\nabla^\omega}(\gamma)}}^{-1}s)={{hol_{\nabla^\omega}(\gamma)}}^{-1}\boldsymbol\alpha \ ,
\end{eqnarray*}and
\begin{eqnarray*}
({{hol_{\nabla^\omega}(\gamma)}}^{-1}-1)\boldsymbol\alpha &=& \boldsymbol{\phi_{2\pi}}^*\boldsymbol\alpha-\boldsymbol{\phi_0}^*\boldsymbol\alpha=\integral{0}{2\pi}{\dev{}{t}(\boldsymbol{\phi_t}^*\boldsymbol\alpha)}{t} \\ &=&\integral{0}{2\pi}{\dev{}{s}\boldsymbol{\phi_{t+s}}^*\boldsymbol\alpha\Big{|}_{s=0}}{t}=\integral{0}{2\pi}{\dev{}{s}\boldsymbol{\phi_s}^*(\boldsymbol{\phi_t}^*\boldsymbol\alpha)\Big{|}_{s=0}}{t} \\ &=&\integral{0}{2\pi}{\pounds_X^{\nabla}(\boldsymbol{\phi_t}^*\boldsymbol\alpha)}{t} \\ &=&\integral{0}{2\pi}{(\boldsymbol{i}_X\circ\ud^\nabla+\ud^\nabla\circ\boldsymbol{i}_X)(\boldsymbol{\phi_t}^*\boldsymbol\alpha)}{t} \\
&=&\boldsymbol{i}_X\left(\integral{0}{2\pi}{\ud^\nabla(\boldsymbol{\phi_t}^*\boldsymbol\alpha)}{t}\right)+\ud^\nabla\circ\boldsymbol{i}_X\left(\integral{0}{2\pi}{\boldsymbol{\phi_t}^*\boldsymbol\alpha}{t}\right) . 
\end{eqnarray*}Using that the pullback commutes with the derivative (proposition \ref{pullcommutes}), one gets from the last equation
\begin{equation}\label{homotopy2}
({{hol_{\nabla^\omega}(\gamma)}}^{-1}-1)\boldsymbol\alpha=\boldsymbol{i}_X\left(\integral{0}{2\pi}{\boldsymbol{\phi_t}^*(\ud^\nabla\boldsymbol\alpha)}{t}\right)+\ud^\nabla\circ\boldsymbol{i}_X\left(\integral{0}{2\pi}{\boldsymbol{\phi_t}^*\boldsymbol\alpha}{t}\right) \ ,
\end{equation}which resembles the equation satisfied by a homotopy operator.\par

\begin{prop}\label{almosthomo} The expression $\boldsymbol{J}_X(\boldsymbol\alpha)=\boldsymbol{i}_X\left(\integral{0}{2\pi}{\boldsymbol{\phi_t}^*\boldsymbol\alpha}{t}\right)$ defines a degree $-1$ derivation on ${S_P}^\bullet(L)$.
\end{prop}
\textit{Proof:} Propositions \ref{degreeminus} and \ref{pullcommutes} imply that $\boldsymbol{J}_X$ is a derivation, and the degree comes from the fact that $\boldsymbol{i}_X$ has degree $-1$.\hfill $\blacksquare$\par\vspace{0.5em}

The equation (\ref{homotopy2}) implies that $\boldsymbol{J}_X$ satisfies
\begin{equation}\label{homotopy}
({{hol_{\nabla^\omega}(\gamma)}}^{-1}-1)\boldsymbol\alpha=\boldsymbol{J}_X(\ud^\nabla\boldsymbol\alpha)+\ud^\nabla\boldsymbol{J}_X(\boldsymbol\alpha) \ ,
\end{equation}for any $\boldsymbol\alpha\in S_P^k(L)$ if $k\geq 1$, whilst for $k=0$ it becomes
\begin{equation}\label{homotopy1} 
({{hol_{\nabla^\omega}(\gamma)}}^{-1}-1)\boldsymbol\alpha=\boldsymbol{J}_X(\ud^\nabla\boldsymbol\alpha) \ ,
\end{equation}since $S^{-1}_P(L)$ is empty and $\boldsymbol{J}_X$ has degree $-1$.

\begin{prop}\label{constholonomy} $\ud^\nabla(({{hol_{\nabla^\omega}(\gamma)}}^{-1}-1)\boldsymbol\alpha)=({{hol_{\nabla^\omega}(\gamma)}}^{-1}-1)\ud^\nabla\boldsymbol\alpha$ for any $\boldsymbol\alpha\in S^k_P(L)$; hence, ${hol_{\nabla^\omega}(\gamma)}$ is constant along $P$.
\end{prop}
\textit{Proof:} It is a direct consequence of equation (\ref{homotopy}): 
\begin{equation}
\ud^\nabla(({{hol_{\nabla^\omega}(\gamma)}}^{-1}-1)\boldsymbol\alpha)=\ud^\nabla[\boldsymbol{J}_X(\ud^\nabla\boldsymbol\alpha)+\ud^\nabla\boldsymbol{J}_X(\boldsymbol\alpha)]=\ud^\nabla\boldsymbol{J}_X(\ud^\nabla\boldsymbol\alpha) \ ,
\end{equation} 
\begin{equation}
({{hol_{\nabla^\omega}(\gamma)}}^{-1}-1)\ud^\nabla\boldsymbol\alpha=\boldsymbol{J}_X(\ud^\nabla\circ\ud^\nabla\boldsymbol\alpha)+\ud^\nabla\boldsymbol{J}_X(\ud^\nabla\boldsymbol\alpha)=\ud^\nabla\boldsymbol{J}_X(\ud^\nabla\boldsymbol\alpha) \ .
\end{equation} \hfill $\blacksquare$\par\vspace{0.5em}

\begin{lem}\label{monodromy}Let $X$ be the generator of a symplectic $S^1$-action; then, 
\begin{equation}\label{holonomyformula}
{hol_{\nabla^\omega}(\gamma)}=c([\gamma])\cdot\e^{i2\pi\theta(X)} \ ,
\end{equation}where $\theta$ is a particular invariant potential $1$-form for $\omega$ in a neighbourhood of $\gamma$ and $c\in\mathrm{Hom}(\pi_1(M);U(1))$. 
\end{lem}
\textit{Proof:} Weinstein's theorem for isotropic embeddings \cite{Wein} asserts that in a neighbourhood $N$ of an orbit the symplectic form is exact, $\omega=\ud\theta$ ---the potential $1$-form can be chosen to be invariant by averaging it with respect to the flow of $X$. Let $s\in\Gamma(\CC\times N)$ be a unitary section of the trivial bundle given by lemma \ref{flatbundle} which has $\theta$ as the potential $1$-form for the nonflat part of $\nabla^\omega$.

Cartan's magic formula and the invariance of $\theta$ give: 
\begin{equation}
0=\pounds_X(\theta)=\imath_X\ud\theta+\ud(\imath_X\theta) \ \Rightarrow \ \imath_X\omega=-\ud\theta(X) \ ;
\end{equation}wherefore, near $\gamma$, the action is hamiltonian, and $\theta(X)$ is its hamiltonian function. In particular, since $\gamma$ is an integral curve of the hamiltonian flow, $\theta(\dot{\gamma}(t))$ is constant.

The parallel transport of a section $r=fs$, with $f\in C^\infty(N)$, around the loop $\gamma$ is given by solving the equation $\nabla^\theta_{\dot{\gamma}}(r)=0$, which is equivalent to
\begin{equation}\label{paralelo}
\dev{}{t}f\circ\gamma(t)=i\theta(\dot{\gamma}(t))f\circ\gamma(t) \ .
\end{equation}Indeed,
\begin{equation}
0=\nabla^\theta_{\dot{\gamma}}(r)=\dot{\gamma}(f)s-if\theta(\dot{\gamma})s
\end{equation}and $\dot{\gamma}(f)\big{|}_{\gamma(t)}=\devv{}{t}f\circ\gamma(t)$; thus, 
\begin{equation}
[\nabla^\theta_{\dot{\gamma}}(r)](\gamma(t))=\left(\dev{}{t}f\circ\gamma(t)\right)s\circ\gamma(t)-if\circ\gamma(t)\theta(\dot{\gamma}(t))s\circ\gamma(t) \ .
\end{equation}

Because $\theta(\dot{\gamma}(t))$ is constant, $f\circ\gamma(t)=\e^{it\theta(X)}$ is a solution to equation (\ref{paralelo}); therefore, if $c([\gamma])$ is the holonomy of the flat hermitian line bundle $(E,\nabla^0)$, then 
\begin{equation}
{hol_{\nabla^\omega}(\gamma)}=c([\gamma])\cdot\e^{i2\pi\theta(X)} \ .
\end{equation}\hfill $\blacksquare$\par\vspace{0.5em}

\begin{rema}
The set of points in $M$ fixed by the $S^1$-action, $M^{S^1}$, is a closed submanifold of $M$, and lemma \ref{monodromy} implies that critical points of ${hol_{\nabla^\omega}(\gamma)}$ are fixed points of the $S^1$-action. If the action is nontrivial, the complement of $M^{S^1}$ is an open dense subset of $M$, and on this subset the action is free; thus, $\{{hol_{\nabla^\omega}(\gamma)}=1\}\cap(M-M^{S^1})$ is a codimension $1$ submanifold.  This means that the set of points with trivial holonomy is a stratified submanifold, and its top dimensional strata have codimension $1$. 
\end{rema}

\begin{prop}\label{rawaction} Supposing that $(M,\omega)$ admits a nontrivial symplectic $S^1$-action whose generator belongs to $P$, flat sections of $L$ vanish.
\end{prop}
\textit{Proof:} Let $s\in\Gamma(L)$ be a flat section, $\nabla s=0$. By equation (\ref{homotopy1}),  
\begin{equation*}
({{hol_{\nabla^\omega}(\gamma)}}^{-1}-1)s=0
\end{equation*}and the flat section vanishes on the dense set where ${hol_{\nabla^\omega}(\gamma)}\neq 1$. Consequently, there are no nonzero flat sections.\hfill $\blacksquare$\par\vspace{0.5em}

\begin{lem}\label{rawmagic}Supposing that $(M,\omega)$ admits a nontrivial symplectic $S^1$-action whose generator belongs to $P$, a form $\boldsymbol\alpha\in S^k_P(L)$ vanishes where ${hol_{\nabla^\omega}(\gamma)}=1$ if and only if there exists a $\boldsymbol\beta\in S^k_P(L)$ such that $\boldsymbol\alpha=({{hol_{\nabla^\omega}(\gamma)}}^{-1}-1)\boldsymbol\beta$.
\end{lem}
\textit{Proof:}  
If $\boldsymbol\alpha=({{hol_{\nabla^\omega}(\gamma)}}^{-1}-1)\boldsymbol\beta$ it is obvious that $\boldsymbol\alpha$ vanishes where ${hol_{\nabla^\omega}(\gamma)}=1$. If the converse holds for functions on $M$, in any trivialising neighbourhood $A$ with unitary section $s$ and coordinates $(z_1,\dots,z_{2n})$, the form $\boldsymbol\alpha$ can be expressed by 
\begin{equation}
\boldsymbol\alpha=\left[\somatorio{j_1,\dots,j_k=1}{2n}\alpha_{j_1,\dots,j_k}(z_1,\dots,z_{2n})\ud z_{j_1}\wedge\cdots\wedge\ud z_{j_k}\right]\otimes s \ .
\end{equation}Furthermore, $\boldsymbol\alpha=0$ at $\{{hol_{\nabla^\omega}(\gamma)}=1\}$ if and only if all the functions $\alpha_{j_1,\dots,j_k}$ vanish on $A\cap\{{hol_{\nabla^\omega}(\gamma)}=1\}$. Hence, there exist functions $\beta_{j_1,\dots,j_k}$ such that $\alpha_{j_1,\dots,j_k}=({{hol_{\nabla^\omega}(\gamma)}}^{-1}-1)\beta_{j_1,\dots,j_k}$. The manifold $M$ can be covered by trivialising neighbourhoods, and the local functions $\beta_{j_1,\dots,j_k}$ piece together to give the desired $\boldsymbol\beta\in S^k_P(L)$. 

Therefore, given $f\in C^{\infty}(A)$ satisfying $f|_{A\cap\{{hol_{\nabla^\omega}(\gamma)}=1\}}=0$ one must construct a $g\in C^{\infty}(A)$ such that $f=({hol_{\nabla^\omega}(\gamma)}^{-1}-1)g$.\par

For points where $1$ is a regular value of ${hol_{\nabla^\omega}(\gamma)}$, theorem 4 in \cite{JHR} proves that this expression holds for functions. On the other hand, lemma \ref{monodromy} implies that locally ${hol_{\nabla^\omega}(\gamma)}=\e^{2\pi ih}$ for some function $h$. 

Let $A$ be a neighbourhood of a critical point of $\{{hol_{\nabla^\omega}(\gamma)}=1\}$. By shrinking $A$, and possibly changing $h$ by a constant, one can assume that only $0$, and no other integer, satisfies $A\cap h^{-1}(\{0\})\neq\emptyset$. With the aid of the flow $\varphi_t$ of the vector field $-hZ$, where in local coordinates around the critical point $Z=h\somatorio{j=1}{2n} z_j\depp{}{z_j}$, one can define a function $g\in C^{\infty}(A)$: 
\begin{equation}\label{eqrawmagic}
g=\fraction{\integral{0}{\infty}{Z(f\circ\varphi_t)}{t}}{2\pi i\integral{0}{1}{\e^{-2\pi ith}}{t}} \ .
\end{equation}In fact, for $h=0$
\begin{equation}
\integral{0}{1}{\e^{-2\pi ith}}{t}=1 \ ,
\end{equation}and for $h\neq 0$
\begin{equation}
\integral{0}{1}{\e^{-2\pi ith}}{t}=\fraction{{hol_{\nabla^\omega}(\gamma)}^{-1}-1}{-2\pi ih} \ .
\end{equation}Thus, the denominator in expression (\ref{eqrawmagic}) never vanishes, whilst 
\begin{eqnarray}
g&=&\fraction{\integral{0}{\infty}{Z(f\circ\varphi_t)}{t}}{2\pi i({{hol_{\nabla^\omega}(\gamma)}}^{-1}-1)/(-2\pi ih)}=\fraction{-\integral{0}{\infty}{hZ(f\circ\varphi_t)}{t}}{{hol_{\nabla^\omega}(\gamma)}^{-1}-1} \nonumber \\ &=&\fraction{\integral{\infty}{0}{\dev{}{t}f\circ\varphi_t}{t}}{{hol_{\nabla^\omega}(\gamma)}^{-1}-1}=\fraction{f-\limite{t}{\infty}f\circ\varphi_t}{{hol_{\nabla^\omega}(\gamma)}^{-1}-1} \ .
\end{eqnarray}\par 

For any point $p\in A$ the limit $\limite{t}{\infty}\varphi_t(p)$ is the critical point (which, in particular, has trivial holonomy) and $f|_{A\cap\{{hol_{\nabla^\omega}(\gamma)}=1\}}=0$; consequently, $f=({{hol_{\nabla^\omega}(\gamma)}}^{-1}-1)g$ on $A\cap\{{hol_{\nabla^\omega}(\gamma)}\neq 1\}$. By continuity of $f,g$ and density of $\{{hol_{\nabla^\omega}(\gamma)}\neq 1\}$, this must be true over all $A$.\hfill $\blacksquare$\par\vspace{0.5em}

The next proposition is a key tool to prove that the Kostant complex is a fine resolution, theorem \ref{fineresolution1}, when the (singular) polarisation comes from an almost or locally toric structure.

\begin{prop}\label{rawpotential} Supposing that $(M,\omega)$ admits a nontrivial symplectic $S^1$-action whose generator belongs to $P$, let $\boldsymbol\alpha\in S^k_P(L)$ be closed, $\ud^\nabla\boldsymbol\alpha=0$, and $k\neq 0$. 
\begin{itemize}
\item The form $\boldsymbol\alpha$ is exact everywhere ${hol_{\nabla^\omega}(\gamma)}\neq 1$. It is also globally exact if $\boldsymbol{J}_X(\boldsymbol\alpha)=0$ where ${hol_{\nabla^\omega}(\gamma)}=1$. 
\item When $\{{hol_{\nabla^\omega}(\gamma)}=1\}$ is a (not necessarily connected) submanifold, $\boldsymbol\alpha$ is exact on $M$ if and only if $\boldsymbol{J}_X(\boldsymbol\alpha)\big{|}_{T\{{hol_{\nabla^\omega}(\gamma)}=1\}}$ is exact.
\end{itemize}
\end{prop}
\textit{Proof:} At points satisfying ${hol_{\nabla^\omega}(\gamma)}\neq 1$ a $(k-1)$-form $\boldsymbol\beta$ is well defined by 
\begin{equation}\label{potentialbeta}
\boldsymbol\beta=\frac{\boldsymbol{J}_X(\boldsymbol\alpha)}{{{hol_{\nabla^\omega}(\gamma)}}^{-1}-1} \ .
\end{equation}Proposition \ref{constholonomy} and equation (\ref{homotopy}), together with the hypothesis of $\boldsymbol\alpha$ being closed, imply that $\ud^\nabla\boldsymbol\beta=\boldsymbol\alpha$. In other words, $\boldsymbol{J}_X/({{hol_{\nabla^\omega}(\gamma)}}^{-1}-1)$ is a homotopy operator where ${hol_{\nabla^\omega}(\gamma)}\neq 1$.\par

For $\boldsymbol{J}_X(\boldsymbol\alpha)=0$ at $\{{hol_{\nabla^\omega}(\gamma)}=1\}$, lemma \ref{rawmagic} gives a $\boldsymbol\sigma\in S^{k-1}_P(L)$ such that $\boldsymbol{J}_X(\boldsymbol\alpha)=({{hol_{\nabla^\omega}(\gamma)}}^{-1}-1)\boldsymbol\sigma$; therefore,  $\boldsymbol\beta$ is well defined by the expression (\ref{potentialbeta}).\par

Assuming that $\{{hol_{\nabla^\omega}(\gamma)}=1\}$ is a submanifold, one consequence of proposition \ref{constholonomy} (as it was observed in \cite{JHR}) is that the polarisation is tangent to it, and all definitions make sense with $M$ replaced by $\{{hol_{\nabla^\omega}(\gamma)}=1\}$.\par 

If $\boldsymbol\alpha=\ud^\nabla\boldsymbol\beta$, by applying equation (\ref{homotopy}),
\begin{equation}
\boldsymbol{J}_X(\boldsymbol\alpha)=\boldsymbol{J}_X\circ\ud^\nabla\boldsymbol\beta=({{hol_{\nabla^\omega}(\gamma)}}^{-1}-1)\ud^\nabla\boldsymbol\beta-\ud^\nabla\circ\boldsymbol{J}_X(\boldsymbol\beta) \ ; 
\end{equation}and $\boldsymbol{J}_X(\boldsymbol\alpha)\big{|}_{T\{{hol_{\nabla^\omega}(\gamma)}=1\}}$ is exact.\par 

Conversely, if $\boldsymbol{J}_X(\boldsymbol\alpha)\big{|}_{T\{{hol_{\nabla^\omega}(\gamma)}=1\}}=\ud^\nabla\big{|}_{T\{{hol_{\nabla^\omega}(\gamma)}=1\}}\boldsymbol\zeta$, taking an extension $\boldsymbol\eta\in S^{k-2}_P(L)$ of $\boldsymbol\zeta$, the formula $(\boldsymbol{J}_X(\boldsymbol\alpha)-\ud^\nabla\boldsymbol\eta)\big{|}_{T\{{hol_{\nabla^\omega}(\gamma)}=1\}}=0$ holds and lemma \ref{rawmagic} provides a $\boldsymbol\beta\in S^{k-1}_P(L)$ such that $\boldsymbol{J}_X(\boldsymbol\alpha)-\ud^\nabla\boldsymbol\eta=({{hol_{\nabla^\omega}(\gamma)}}^{-1}-1)\boldsymbol\beta$. Proposition \ref{constholonomy} implies that
\begin{equation}
\ud^\nabla\circ\boldsymbol{J}_X(\boldsymbol\alpha)=\ud^\nabla(({{hol_{\nabla^\omega}(\gamma)}}^{-1}-1)\boldsymbol\beta)=({{hol_{\nabla^\omega}(\gamma)}}^{-1}-1)\ud^\nabla\boldsymbol\beta \ , 
\end{equation}but equation (\ref{homotopy}) reads
\begin{equation}
\ud^\nabla\circ\boldsymbol{J}_X(\boldsymbol\alpha)=({{hol_{\nabla^\omega}(\gamma)}}^{-1}-1)\boldsymbol\alpha \ ;
\end{equation}
thus, $\ud^\nabla\boldsymbol\beta=\boldsymbol\alpha$ holds where ${hol_{\nabla^\omega}(\gamma)}\neq 1$. Since $\ud^\nabla\boldsymbol\beta$ is everywhere defined and $\{{hol_{\nabla^\omega}(\gamma)}\neq 1\}$ is a dense set, $\boldsymbol\alpha$ must be exact.\hfill $\blacksquare$\par\vspace{0.5em}     


\section{The Bohr-Sommerfeld condition}

\hspace{1.5em}The following definition plays a very important role in the computation of the cohomology groups appearing in geometric quantisation:

\begin{defi}\label{B-Sleaf}A leaf $\ell$ of $\mathcal{P}$ is a Bohr-Sommerfeld leaf if there exists a nonzero section $s:\ell\to L$ such that $\nabla^\omega_X s=0$ for any vector field $X$ of $P$ restricted to $\ell$. 
\end{defi}

\begin{prop}\label{bsleaf}A leaf $\ell$ of $\mathcal{P}$ is a Bohr-Sommerfeld leaf if and only if the holonomy is trivial, ${hol_{\nabla^\omega}(\gamma)}=1$, for any loop $\gamma$ on a connected component of $\ell$.\end{prop}
\textit{Proof:} In a Bohr-Sommerfeld leaf $\ell$ the nonzero section $s$ can be used to define a potential $1$-form $\Theta$ of the connexion on the whole leaf. The potential $1$-form vanishes on $\ell$, since $0=\nabla^\omega s\big{|}_{T\ell}=-i\Theta\big{|}_{T\ell}\otimes s$. Thus, if $\gamma$ is a loop on $\ell$, by lemma \ref{monodromy}, ${hol_{\nabla^\omega}(\gamma)}=1$.\par

Now, supposing that ${hol_{\nabla^\omega}(\gamma)}=1$ for any loop on a connected component of a leaf $\ell$ of $\mathcal{P}$, for any point $p\in\ell$ and a nonzero $s_p\in L_p$ (the fibre of $L$ over $p$) it is possible to define a nonzero section $s$ over $\ell$ by parallel transport: i.e. $s(q)=\Pi_{\gamma_1}(s_p)$, where $\gamma_1$ is any curve connecting $p$ and $q\in\ell$. The section is well-defined because if $\gamma_2$ is another curve connecting $p$ and $q$, and $\gamma$ the loop formed by composing $\gamma_2$ and $\gamma_1^{-1}$, 
\begin{eqnarray}
s(q) &=& {hol_{\nabla^\omega}(\gamma)}s(q)=\Pi_{\gamma}(s(q))=\Pi_{\gamma_2}\circ[\Pi_{\gamma_1}]^{-1}(s(q)) \nonumber \\ &=& \Pi_{\gamma_2}\circ[\Pi_{\gamma_1}]^{-1}\circ\Pi_{\gamma_1}(s_p)=\Pi_{\gamma_2}(s_p) \ .
\end{eqnarray}The parallel transport respects the hermitian product, and this guarantees that the section defined in this way is nonzero. \hfill $\blacksquare$\par\vspace{0.5em}

For the example \ref{cilindro}, Bohr-Sommerfeld leaves are circles of integral height. Similarly for the complex plane, endowed with the canonical symplectic structure and a polarisation induced by a circle action.

\begin{exam}\label{complexplane} Let $M=\CC$ with coordinates $(x,y)$ and Darboux form $\omega=\ud x\wedge\ud y$, $L=\CC\times\CC$ the trivial bundle with connexion 1-form $\Theta=\frac{1}{2}(x\ud y-y\ud x)$, with respect to the unitary section $\e^{i(x^2+y^2)}$, and $P=\left\langle -y\depp{}{x}+x\depp{}{y}\right\rangle_{C^\infty(\CC)}$. 

Lemma \ref{monodromy} implies that the holonomy for a closed curve $\gamma$ inside a leaf of $\mathcal{P}$ is 
\begin{equation}
[hol_{\nabla^\omega}(\gamma)](x,y)=\e^{i 2\pi \frac{(x^2+y^2)}{2}} \ .
\end{equation}Therefore, by proposition \ref{bsleaf}, flat sections are only well-defined for the set of points with $\frac{(x^2+y^2)}{2}\in\CZ$.\hfill $\Diamond$\par\vspace{0.5em}
\end{exam}

There is a stronger characterisation for the Bohr-Sommerfeld leaves in the case of integrable systems: this is an application of lemma \ref{monodromy}.

The preimage of a point in $\CR^n$ by the moment map $F:M\to\CR^n$ of an integrable system on a symplectic manifold $(M,\omega)$ of dimension $2n$ is a fibre of the fibration associated with the foliation induced by the distribution generated by the hamiltonian vector fields of the moment map: i.e. the preimage of a point by a moment map is a union of isotropic leaves. The preimage of the origin of $\CR^n$ is called the zero fibre, and fibres which are a union of Bohr-Sommerfeld leaves are called Bohr-Sommerfeld fibres.

\begin{teo}\label{integerlattice2}Under the assumption that the zero fibre is Bohr-Sommerfeld, the image of Bohr-Sommerfeld fibres by a moment map is contained in $\CR^{n-k}\times\CZ^k$; $k$ being the number of linearly independent hamiltonian $S^1$-actions generated by the moment map. 
\end{teo}
\textit{Proof:} This was proved by Guillemin and Sternberg in \cite{GuSt} when the fibres are Liouville tori ---their proof holds for lagrangian fibrations with compact connected fibres over simply connected basis. Lemma \ref{monodromy} and proposition \ref{bsleaf} imply that over a Bohr-Sommerfeld fibre each component of the moment map generating a $S^1$-action takes an integral value, depending only on the fibre.\hfill $\blacksquare$\par\vspace{0.5em} 

\begin{exam} For toric manifolds the Bohr-Sommerfeld leaves are the inverse image by the moment map of integer lattice points in the polytope, with regular ones inside the polytope. \hfill $\Diamond$\par\vspace{0.5em}
\end{exam}


\section{Applications I: local and semilocal computations}


\subsection{The cylinder: polarisation by circles}\label{Mcilindro}

\hspace{1.5em}Recalling example \ref{cilindro}: $(M=\CR\times S^1,\omega=\ud x\wedge\ud y)$, $L$ the trivial bundle with connexion 1-form $\Theta=x\ud y$, with respect to the unitary section $\e^{ix}$, and $P=\left\langle\depp{}{y}\right\rangle_{C^\infty(\CR\times S^1)}$.\par

The hamiltonian vector field $X=\depp{}{y}$ generates a $S^1$-action, and the holonomy of its orbits is given by $[{hol_{\nabla^\omega}(\gamma)}](x,y)=\e^{i2\pi x}$ (lemma \ref{monodromy}, for $S^1\cong\CR/2\pi\CZ$). Proposition \ref{rawaction} holds, and one gets $\check{H}^0(M;\mathcal{J})=\{0\}$ by applying theorem \ref{fineresolution1}. Furthermore, proposition \ref{rawpotential} and theorem \ref{fineresolution1} can be applied, implying $\check{H}^l(V;\mathcal{J}\big{|}_V)=\{0\}$, for $l\geq 1$, for each neighbourhood $V=(a,b)\times S^1$ that does not contain a Bohr-Sommerfeld leaf.

Let $\ell_k$ be the inverse image by the height function of the point $x=k\in\CZ$. Wherefore, $\ell_k\cong S^1$ is a Bohr-Sommerfeld leaf and $\{{hol_{\nabla^\omega}(\gamma)}=1 \}=\displaystyle\bigcup_{k\in\CZ}\ell_k$.\par

It is possible\footnote{This construction is due to Rawnsley \cite{JHR}.} to define a linear map $\Psi:S^1_P(L)\to\displaystyle\bigoplus_{k\in\CZ}\Gamma(L|_{\ell_k})$ by:
\begin{equation}
\Psi(\boldsymbol\alpha)=\oplus_{k\in\CZ}\boldsymbol{J}_X(\boldsymbol\alpha)\Big{|}_{\ell_k} \ .
\end{equation}Because the dimension of $M$ is $2$, $S^l_P(L)=\{0\}$ for $l\geq 2$, and, for any $\boldsymbol\alpha\in S^1_P(L)$, equation (\ref{homotopy}) reads  
\begin{equation}
\nabla\circ\boldsymbol{J}_X(\boldsymbol\alpha)=({{hol_{\nabla^\omega}(\gamma)}}^{-1}-1)\boldsymbol\alpha \ \Rightarrow \ \nabla\Psi(\boldsymbol\alpha)=0 \ .
\end{equation}Thus, the image of $\Psi$ is contained in the set of flat sections over Bohr-Sommerfeld leaves.\par  

Conversely, given $\oplus_{k\in\CZ}s_k\in\displaystyle\bigoplus_{k\in\CZ}\Gamma(L|_{\ell_k})$, where $s_k$ are flat sections ($\nabla s_k=0$), there exists $s\in\Gamma(L)$ such that $s|_{\ell_k}=s_k$ for each $k\in\CZ$: due the closedness of $\displaystyle\bigcup_{k\in\CZ}\ell_k$. Lemma \ref{rawmagic} implies the existence of an $\boldsymbol\alpha\in S^1_P(L)$ satisfying 
\begin{eqnarray}
\nabla s&=&({{hol_{\nabla^\omega}(\gamma)}}^{-1}-1)\boldsymbol\alpha \ \Rightarrow \nonumber \\
({{hol_{\nabla^\omega}(\gamma)}}^{-1}-1)\boldsymbol{J}_X(\boldsymbol\alpha)&=&\boldsymbol{J}_X(\nabla s)=({{hol_{\nabla^\omega}(\gamma)}}^{-1}-1)s \ .
\end{eqnarray}By density and continuity, $\boldsymbol{J}_X(\boldsymbol\alpha)=s$; hence, the image of $\Psi$ is the set of flat sections over Bohr-Sommerfeld leaves.\par    

Proposition \ref{rawpotential} asserts that $\mathrm{ker}\Psi=\nabla(\Gamma(L))$, and the first isomorphism theorem 
\begin{equation}
\begin{array}{c c c}
S^1_P(L) & \longrightarrow  & \Psi(S^1_P(L)) \\
\big{\downarrow} & \swarrow &  \\
S^1_P(L)/\mathrm{ker}\Psi &  & 
\end{array} 
\end{equation}implies that $\check{H}^1(M;\mathcal{J})\cong\displaystyle\bigoplus_{k\in\CZ}\CC$: the ring of flat sections over $\ell_k$ is isomorphic to $\CC$ (see example \ref{cilindro}).

\begin{prop}\label{Qcilindro}The quantisation of a cylinder polarised by circles is $\CC^{b_s}$, where $b_s$ is the number of Bohr-Sommerfeld leaves.   
\end{prop}


\subsection{The complex plane: polarisation by circles}\label{Mdisk}

\hspace{1.5em}Let $(M=\CC,\omega=\ud x\wedge\ud y)$ and $F:M\to\CR$ be a nondegenerate integrable system of elliptic type, i.e. $F(x,y)=x^2+y^2$. For this case, the polarisation is $P=\left\langle-y\depp{}{x}+x\depp{}{y}\right\rangle_{c^\infty(\CC)}$ and the hamiltonian vector field $X=-y\depp{}{x}+x\depp{}{y}$ is a generator of a $S^1$-action ---this is example \ref{complexplane}, again.\par

As in the previous cases, $(M,\omega)$ is an exact symplectic manifold and the trivial line bundle is a prequantum line bundle for it: $L=\CC\times\CC$ with connexion 1-form $\Theta=\frac{1}{2}(x\ud y-y\ud x)$, with respect to the unitary section $\e^{i(x^2+y^2)}$.\par 

\begin{prop}\label{poincaeli1} At a sufficiently small neighbourhood of the origin of $\CC$, the cohomology groups $H^k({S_P}^\bullet(L))$ vanish for $k\geq 0$ when $\mathcal{P}$ is generated by $-y\depp{}{x}+x\depp{}{y}$.
\end{prop}
\textit{Proof:} Proposition \ref{rawaction} can be applied to prove that $H^0({S_P}^\bullet(L))=\{0\}$.

Lemma \ref{monodromy} asserts that $[{hol_{\nabla^\omega}(\gamma)}](x,y)=\e^{i\pi (x^2+y^2)}$; as a result, the set \linebreak $\{{hol_{\nabla^\omega}(\gamma)}=1\}$ is the union of the origin and concentric circles with $R^2/2$ integer ($R$ being the radius), and since the origin is a fixed point, the operator $\boldsymbol{J}_X$ is the null operator when restricted to the origin. Hence, proposition \ref{rawpotential}, applied for each contractible neighbourhood of the origin that does not contain any other Bohr-Sommerfeld leaf, implies that elliptic singularities give no contribution to quantisation\footnote{This was first proved in \cite{Ha} using different techniques.}, $H^k({S_P}^\bullet(L))=\{0\}$ for $k\geq 1$. \hfill$\blacksquare$\par\vspace{0.5em}

Rephrasing the previous proposition, theorem \ref{fineresolution1} holds in this particular setting.

\begin{prop}\label{Qelliptic}The quantisation of an open disk polarised by circles is $\CC^{b_s}$, where $b_s$ is the number of nonsingular Bohr-Sommerfeld leaves.  
\end{prop}
\textit{Proof:} Let $M$ be the complex plane; for an arbitrary open disk the same argument works. $M$ can be divided up into an open disk $V$ of radius $b<1$ centred at the origin, and an annulus $W$ centred at the origin with small radius $a\in (0,b)$ and an infinite big radius: $M=V\cup W$ and $V\cap W$ is an annulus with small radius $a$ and big radius $b$. Proposition \ref{poincaeli1} implies that $\check{H}^k(V;\mathcal{J}\big{|}_V)=\{0\}$ for all $k$, and proposition \ref{Qcilindro} gives $\check{H}^k(V\cap W;\mathcal{J}\big{|}_{V\cap W})=\{0\}$ for all $k$ as well, since $V\cap W\cong (a,b)\times S^1$ (polarised by circles). Thus, the Mayer-Vietoris argument works and $\check{H}^k(M;\mathcal{J})\cong\check{H}^k(W;\mathcal{J}\big{|}_W)$. It happens that $W\cong (a,\infty)\times S^1$ (polarised by circles), and proposition \ref{Qcilindro} concludes the proof. \hfill$\blacksquare$\par\vspace{0.5em}


\subsection{Symplectic vector spaces: linear polarisation}

\hspace{1.5em}All symplectic vector spaces polarised by lagrangian hyperplanes are equivalent to this particular case, this is the local model for polarised symplectic manifolds when the real polarisation has no singularities: 
\begin{equation*}
(\CC^n,\somatorio{j=1}{n}\ud x_j\wedge\ud y_j) \ \text{and polarisation} \  \left\langle\depp{}{y_1},\dots,\depp{}{y_n}\right\rangle_{C^\infty(\CC^n)} \ .
\end{equation*}

Since $(\CC^n,\somatorio{j=1}{n}\ud x_j\wedge\ud y_j)$ is an exact symplectic manifold, the trivial line bundle is a prequantum line bundle for it: $\CC\times\CC^n$ with connexion 1-form $\somatorio{j=1}{n} x_j\ud y_j$, with respect to the unitary section $\exp\left( i\somatorio{j=1}{n} x_j\right)$.

The solutions of the flat equation are complex-valued functions of the type 
\begin{equation}
h(x_1,\dots,x_n)\exp\left( i\somatorio{j=1}{n} x_jy_j\right) \ .
\end{equation}Therefore, 
\begin{equation}
\check{H}^0(\CC^n;\mathcal{J})\cong C^\infty(\CR^n) \ .
\end{equation}  

\begin{lem}\label{locflat} There always exists a local unitary flat section on each point of $M$ for a given regular polarisation $P$.
\end{lem}

As a consequence of the existence of unitary flat sections, elements of $\mathcal{S}^k_P(L)$ which are closed can be interpreted as closed polarised $k$-forms taking values on the sheaf $\mathcal{J}$. 

\begin{coro}\label{corofran} Let $\boldsymbol{\Omega^k_P}$ be the sheaf associated to $\Omega^k_P(M)$; then, $\mathcal{S}^k_P(L)\cong\boldsymbol{\Omega^k_P}\otimes\mathcal{J}$ and $\mathrm{ker}(\ud_\nabla)\cong\mathrm{ker}(\ud_P)\otimes\mathcal{J}$.
\end{coro}
\textit{Proof:} By lemma \ref{locflat}, for each point on $M$ there exists a trivialising neighbourhood $V\subset M$ of $L$ with an unitary flat section $s\in\Gamma(L\big{|}_V)$. If $\boldsymbol\alpha\in S^k_P(L)$, it can be locally written as $\boldsymbol\alpha\big{|}_{TV}=\alpha\otimes s$, where $\alpha\in\Omega^k_P(V)$. The condition $\ud^\nabla(\alpha\otimes s)=0$ is, then, equivalent to $\ud_P\alpha=0$, because $\ud^\nabla(\alpha\otimes s)=\ud_P\alpha\otimes s+(-1)^k\alpha\wedge\nabla s$, $s\neq 0$ and $\nabla s=0$.\hfill $\blacksquare$\par\vspace{0.5em}

The Kostant complex \eqref{eq72} is just the foliated complex (associated to the foliation induced by the polarisation) twisted by the sheaf of sections $\mathcal{S}$, and the exactness of the foliated complex implies, by corollary \ref{corofran}, the exactness of the Kostant complex.  

Indeed, this can be made more explicit. Using the local unitary flat section $r=\exp\left( i\somatorio{j=1}{n} x_jy_j\right)$ as basis, if $\alpha\otimes r\in S_P^k(L)$ is closed:
\begin{equation}
0=\ud^\nabla(\alpha\otimes r)=\ud_P\alpha\otimes r+(-1)^k\alpha\wedge\nabla r=\ud_P\alpha\otimes r \ .
\end{equation}Wherefore, Poincaré lemma for regular foliations imply that $H^k({S_P}^\bullet(L))=\{0\}$ for $k\geq 1$.

\begin{teo}\label{poincaregu} At a sufficiently small neighbourhood of any point of $M$, the cohomology groups $H^k({S_P}^\bullet(L))$ vanish for $k\geq 1$ when $\mathcal{P}$ is a subbundle of $TM$.
\end{teo}

From the results of this subsection one concludes that:

\begin{prop}\label{Qplano}The quantisation of the cotangent bundle of $\CR^n$ with linear polarisation is $C^\infty(\CR^n)$.   
\end{prop}

This is by no means an example where the techniques developed in this article are used; howbeit, this result is needed below.


\subsection{Direct product type with a regular component}\label{Kregular}

\hspace{1.5em}The following quantisation problem will be considered now: $N=(-1,1)\times S^1$ endowed with the same structures (symplectic form, polarisation and prequantum line bundle) of the model in subsection \ref{Mcilindro}, and $(M,\omega)$ a prequantisable symplectic manifold with polarisation $P$ and prequantum line bundle $(L,\nabla^\omega)$. The product $N\times M$ admits $\mathscr{P}=\left\langle\depp{}{y}\right\rangle_{C^\infty(N)}\oplus_{C^\infty(N\times M)}P$ as a polarisation for the symplectic form $\ud x\wedge\ud y +\omega$, and also a prequantum line bundle $(\mathscr{L},\bar{\nabla}^{\ud x\wedge\ud y +\omega})$.

The vector field $\depp{}{y}$ generates a hamiltonian $S^1$-action: the holonomy of its orbits is given by ${hol_{\bar{\nabla}^{\ud x\wedge\ud y +\omega}}(\gamma)}=\e^{i2\pi x}$ (lemma \ref{monodromy}). Wherefore, proposition \ref{rawaction} can be used to show that $H^0({S_\mathscr{P}}^\bullet(\mathscr{L}))=\{0\}$.\par

For the other groups one has:

\begin{teo}\label{kunnethregular}Supposing that $(M,\omega)$ is exact, $\omega=\ud\theta$, and both $\mathscr{L}$ and $L$ trivial, the map 
\begin{equation}
\Psi:H^k({S_\mathscr{P}}^\bullet(\mathscr{L}))\to H^{k-1}({S_P}^\bullet(L))
\end{equation}defined by $\Psi([\bar{\boldsymbol\alpha}])=\left[\boldsymbol{J}_{\depp{}{y}}(\bar{\boldsymbol\alpha})\Big{|}_{\{{hol_{\bar{\nabla}^{\ud x\wedge\ud y +\omega}}(\gamma)}=1\}}\right]$ is an isomorphism.
\end{teo}
\textit{Proof:} Since the product symplectic manifold, $(N\times M,\ud x\wedge\ud y +\omega)$, is exact, there exists a unitary section $\bar{s}\in\Gamma(\mathscr{L})$ satisfying $\bar{\nabla}\bar{s}=-i(x\ud y+\theta)\otimes\bar{s}$ (lemma \ref{flatbundle} with $\mathscr{L}$ trivial). Let $\bar{\boldsymbol\alpha}=\bar{\alpha}\otimes\bar{s}\in S^k_\mathscr{P}(\mathscr{L})$ and $\bar{\alpha}=\ud y\wedge\bar{\beta}+\bar{\sigma}$, where $\bar{\beta}=\imath_{\depp{}{y}}\bar{\alpha}$ and $\bar{\sigma}=\bar{\alpha}-\ud y\wedge\bar{\beta}$; thus, $\imath_{\depp{}{y}}\bar{\beta}=\imath_{\depp{}{y}}\bar{\sigma}=0$.\par

\begin{equation}
\boldsymbol{J}_{\depp{}{y}}(\bar{\boldsymbol\alpha})=\integral{0}{2\pi}{\phi^*_t\bar{\beta}\otimes\e^{-itx}\bar{s}}{t}=\integral{0}{2\pi}{\phi^*_t\bar{\beta}\e^{-itx}}{t}\otimes\bar{s} \ \Rightarrow 
\end{equation}
\begin{equation}
\boldsymbol{J}_{\depp{}{y}}(\bar{\boldsymbol\alpha})\Big{|}_{\{{hol_{\bar{\nabla}^{\ud x\wedge\ud y +\omega}}(\gamma)}=1\}}=\eta\otimes\bar{s}|_{x=0} \ ,
\end{equation}where $\eta=\integral{0}{2\pi}{\phi^*_t\bar{\beta}\e^{-itx}}{t}\Big{|}_{\{x=0\}}$. The flow of $\depp{}{y}$ preserves $\eta$; therefore, $\eta\in\Omega^{k-1}_P(M)$.

For closed $\bar{\boldsymbol\alpha}$, 
\begin{equation}
\ud^{\bar{\nabla}}\circ\boldsymbol{J}_{\depp{}{y}}(\bar{\boldsymbol\alpha})=({{hol_{\bar{\nabla}^{\ud x\wedge\ud y +\omega}}(\gamma)}}^{-1}-1)\bar{\boldsymbol\alpha} \ \Rightarrow \ \ud^{\bar{\nabla}}\circ\boldsymbol{J}_{\depp{}{y}}(\bar{\boldsymbol\alpha})\Big{|}_{\{{hol_{\bar{\nabla}^{\ud x\wedge\ud y +\omega}}(\gamma)}=1\}}=0 \ .
\end{equation}

By construction, $\bar{\nabla}_{\depp{}{y}}\bar{s}=-ix\bar{s}$ and $\bar{\nabla}_{\depp{}{y}}\bar{s}|_{x=0}=0$; wherefore, for each point $p\in M$, $\bar{s}|_{x=0}$ is uniquely determined by its value at $(0,0,p)\in N\times M$ by parallel transport along integral curves of $\depp{}{y}$. This means that $\bar{s}|_{x=0}$ identifies itself as a section of $L$: the restriction of $\mathscr{L}$ to $\{(0,0)\}\times M$ is a trivial line bundle over $M$ with a connexion such that its curvature is equal to $\omega$; consequently, it must be isomorphic to $L$.

To summarise it, after some identifications, $\boldsymbol{J}_{\depp{}{y}}(\cdot)\Big{|}_{\{{hol_{\bar{\nabla}^{\ud x\wedge\ud y +\omega}}(\gamma)}=1\}}$ maps closed $k$-forms of ${S_\mathscr{P}}^\bullet(\mathscr{L})$ to closed $(k-1)$-forms of ${S_P}^\bullet(L)$, and proposition \ref{rawpotential} proves that $\Psi$ is injective ---the set $\{{hol_{\bar{\nabla}^{\ud x\wedge\ud y +\omega}}(\gamma)}=1\}$ is equal to $\{0\}\times S^1\times M$.

Now, given $r\in\Gamma(L)$ a unitary section, let $\bar{r}\in\Gamma(\mathscr{L})$ be an extension of the following section defined on $\{x=0\}$: after identifying $\mathscr{L}\big{|}_{\{(0,0)\}\times M}$ with $L$, for each point $p\in M$, the parallel transport of $r(p)$ by the integral curve of $\depp{}{y}$ passing through $p$ defines a section of $\mathscr{L}$ over the set $\{x=0\}$.\par

Due to the inclusion $\Omega^{k-1}_P(M)\subset\Omega^{k-1}_\mathscr{P}(N\times M)$, the expression $\bar{\boldsymbol\zeta}=\zeta\otimes \bar{r}$ defines an element in $S^{k-1}_\mathscr{P}(\mathscr{L})$ for any $[\zeta\otimes r]\in H^{k-1}({S_P}^\bullet(L))$.

The form $\ud^{\bar{\nabla}}\bar{\boldsymbol\zeta}\big{|}_{\{{hol_{\bar{\nabla}^{\ud x\wedge\ud y +\omega}}(\gamma)}=1\}}$ is completely determined by $\ud^{\nabla}(\zeta\otimes r)$, via parallel transport (which commutes with the derivation). And because $\ud^{\nabla}(\zeta\otimes r)=0$, lemma \ref{rawmagic} provides an $\bar{\boldsymbol\alpha}\in S^k_\mathscr{P}(\mathscr{L})$ such that $\ud^{\bar{\nabla}}\bar{\boldsymbol\zeta}=({{hol_{\bar{\nabla}^{\ud x\wedge\ud y +\omega}}(\gamma)}}^{-1}-1)\bar{\boldsymbol\alpha}$. Hence,
\begin{equation}
0=\ud^{\bar{\nabla}}\circ\ud^{\bar{\nabla}}\bar{\boldsymbol\zeta}=({{hol_{\bar{\nabla}^{\ud x\wedge\ud y +\omega}}(\gamma)}}^{-1}-1)\ud^{\bar{\nabla}}\bar{\boldsymbol\alpha} \ ,
\end{equation}implying that $\bar{\boldsymbol\alpha}$ is closed over the dense set $\{x\neq 0\}$, and, by continuity, $\bar{\boldsymbol\alpha}$ is closed.

As consequence of $\imath_{\depp{}{y}}\zeta$ being zero, $\boldsymbol{J}_{\depp{}{y}}(\bar{\boldsymbol\zeta})=0$ and equation (\ref{homotopy}) reads 
\begin{eqnarray}
({{hol_{\bar{\nabla}^{\ud x\wedge\ud y +\omega}}(\gamma)}}^{-1}-1)\bar{\boldsymbol\zeta}&=&\boldsymbol{J}_{\depp{}{y}}\circ\ud^{\bar{\nabla}}\bar{\boldsymbol\zeta}=\boldsymbol{J}_{\depp{}{y}}(({{hol_{\bar{\nabla}^{\ud x\wedge\ud y +\omega}}(\gamma)}}^{-1}-1)\bar{\boldsymbol\alpha}) \nonumber \\ 
&=&({{hol_{\bar{\nabla}^{\ud x\wedge\ud y +\omega}}(\gamma)}}^{-1}-1)\boldsymbol{J}_{\depp{}{y}}(\bar{\boldsymbol\alpha}) \ ,
\end{eqnarray}which implies $\boldsymbol{J}_{\depp{}{y}}(\bar{\boldsymbol\alpha})=\bar{\boldsymbol\zeta}$ where $x\neq 0$; and by density and continuity, it must hold true everywhere. This proves that $\Psi$ is onto.\hfill$\blacksquare$\par\vspace{0.5em}

The theorem still holds if $N$ is replaced by $(a,b)\times S^1$ with $(a,b)\cap\CZ=\{k\}$. For $(a,b)\cap\CZ=\emptyset$, propositions \ref{rawaction} and \ref{rawpotential} assert that all cohomology groups $H^l({S_\mathscr{P}}^\bullet(\mathscr{L}))$ vanish: the quantisation of the product is trivial when there is no Bohr-Sommerfeld leaf. 

By a Mayer-Vietoris argument, similar to the one described above (subsection \ref{Mdisk}), one can compute the product quantisation for $(a,b)\cap\CZ=\{k_1,\dots,k_{b_s}\}$. It suffices to take the cover $\mathcal{A}=\{A_j\}_{j\in\{1,\dots,b_s\}}$, where $A_1=(a,k_1+3/4)\times M$, $A_{b_s}=(k_{b_s}-3/4,b)\times M$ and $A_j=(k_j-3/4,k_j+3/4)\times M$ (supposing $k_1\leq k_2\leq\dots\leq k_{b_s}$).    

\begin{coro}\label{kunnethregular2}Assuming that the Kostant complex is a fine resolution for the sheaf of flat sections of $L$, the quantisation of the product between a cylinder polarised by circles and an arbitrary exact symplectic manifold $M$ is a direct sum of $b_s$ copies of ${\mathcal{Q}(M)}$: where $b_s$ is the number of Bohr-Sommerfeld leaves with respect to the quantisation of the cylinder. 
\end{coro}


\subsection{Direct product type with an elliptic component}

\hspace{1.5em}The quantisation problem to be considered here is: $N=\{(x,y)\in\CC \ ; \ x^2+y^2<1\}$ endowed with the same structures (symplectic form, polarisation and prequantum line bundle) of the model in subsection \ref{Mdisk}, and $(M,\omega)$ a prequantisable symplectic manifold with real polarisation $\mathcal{P}$ and prequantum line bundle $(L,\nabla^\omega)$. The product $N\times M$ admits $\mathscr{P}=\left\langle -y\depp{}{x}+x\depp{}{y}\right\rangle_{C^\infty(N)}\oplus_{C^\infty(N\times M)}P$ as a polarisation for the symplectic form $\ud x\wedge\ud y +\omega$, and also a prequantum line bundle $(\mathscr{L},\bar{\nabla}^{\ud x\wedge\ud y +\omega})$.\par

\begin{lem}\label{poincaeli2} The cohomology groups $H^k({S_\mathscr{P}}^\bullet(\mathscr{L}))$ vanish for $k\geq 0$ in the particular case described in this subsection.
\end{lem}
\textit{Proof:} The group $H^0({S_\mathscr{P}}^\bullet(\mathscr{L}))$ is trivial because $X=-y\depp{}{x}+x\depp{}{y}$ generates a hamiltonian $S^1$-action; wherefore, proposition \ref{rawaction} holds. Whilst for higher order groups, one needs to note that the set $\{{hol_{\bar{\nabla}^{\ud x\wedge\ud y +\omega}}(\gamma)}=1\}$ is equal to $\{(0,0)\}\times M$ (lemma \ref{monodromy} gives ${hol_{\bar{\nabla}^{\ud x\wedge\ud y +\omega}}(\gamma)}=\e^{i\pi(x^2+y^2)}$), and that $(0,0,p)$ are fixed points for any $p\in M$; thus, the operator $\boldsymbol{J}_X$ is the null operator when restricted to $\{(0,0)\}\times M$. Therefore, by applying proposition \ref{rawpotential}, $H^k({S_\mathscr{P}}^\bullet(\mathscr{L}))=\{0\}$ for $k\geq 1$.\hfill$\blacksquare$\par\vspace{0.5em}

By a Mayer-Vietoris argument similar to the ones used in subsections \ref{Mdisk} and \ref{Kregular}, one has:

\begin{prop}\label{kunnethelliptic}Assuming that the Kostant complex is a fine resolution for the sheaf of flat sections of a prequantum line bundle of $(M,\omega)$, the quantisation of the product between an open disk polarised by circles and an arbitrary exact symplectic manifold $M$ is a direct sum of $b_s$ copies of ${\mathcal{Q}(M)}$: where $b_s$ is the number of nonsingular Bohr-Sommerfeld leaves with respect to the quantisation of the open disk.
\end{prop}

The following theorem is the Poincaré lemma for singularities of elliptic type. 

\begin{teo}\label{poincaeli3}
Assuming that $p\in M$ is a nondegenerate critical point of Williamson type $(k_e,k_h,k_f)$, with $k_e\geq 1$ and $k_e+k_h+2k_f\leq n$ (it does not need to be a rank zero critical point, but it has to have an elliptic component), for an integrable system $F:M\to\CR^n$ on a prequantisable symplectic manifold $(M,\omega)$, with polarisation induced by the moment map: the cohomology groups $H^k({S_P}^\bullet(L))$ vanish for $k\geq 0$ in a sufficiently small neighbourhood of $p$. 
\end{teo}
\textit{Proof:} The normal form theorem of Eliasson, Miranda and Zung for nondegenerate singularities of integrable systems (see \cite{BolOsh} and references therein for details) says that the model near a critical point of that type is exactly the one of lemma \ref{poincaeli2}: as long as $k_e\geq 1$, along one of the elliptic directions a neighbourhood of $p$ splits into a product as the model of lemma \ref{poincaeli2}. \hfill$\blacksquare$\par\vspace{0.5em} 


\subsection{Focus-focus singularities}

\hspace{1.5em}Let $F=(f_1,f_2):M\to\CR^2$ be an integrable system on a prequantasible $(M,\omega)$, with a rank zero nondegenerate critical point of Williamson type $(0,0,k_f)$. Near the singular point, $f_2$ generates, via its hamiltonian vector field flow, a hamiltonian $S^1$-action ---Zung \cite{Zu2} demonstrated that this action is defined semilocally, i.e. near a neighbourhood of a focus-focus singular fibre. 

In a small enough neighbourhood $W$ of a singular point of a focus-focus fibre, $\boldsymbol{J}_X$ is the null operator over the points where $\{{hol_{\nabla^\omega}(\gamma)}=1\}$. Indeed, the symplectic local model (provided in \cite{Zu1,Zu2} and attributed to Eliasson, Lerman and Umanskiy, and Vey) is given by, $W\cong\CC^2$ with coordinates $(x_1,x_2,y,_1,y_2)$, $L|_W\cong\CC\times\CC^2$ with connexion 1-form 
\begin{equation}
\Theta=\frac{1}{2}(x_1\ud y_1-y_1\ud x_1+x_2\ud y_2-y_2\ud x_2) \ ,
\end{equation}with respect to the unitary section $s=\e^{i(x_1y_1+x_2y_2+x_1y_2-x_2y_1)}$. 

The integrable system takes the form 
\begin{equation}
F(x_1,x_2,y_1,y_2)=(x_1y_1+x_2y_2, x_1y_2-x_2y_1) \ ,
\end{equation}and, therefore, the polarisation is generated by
\begin{equation}
X_1=-x_1\dep{}{x_1}-x_2\dep{}{x_2}+y_1\dep{}{y_1}+y_2\dep{}{y_2}
\end{equation}and
\begin{equation}
X_2=x_2\dep{}{x_1}-x_1\dep{}{x_2}+y_2\dep{}{y_1}-y_1\dep{}{y_2} \ .
\end{equation}\par

The hamiltonian vector field $X_2$ is the generator of the $S^1$-action. Its periodic flow is given by:
\begin{equation}
\phi_t(x_1,x_2,y_1,y_2)=(x_1\cos t+x_2\sin t,x_2\cos t-x_1\sin t,y_1\cos t+y_2\sin t,y_2\cos t-y_1\sin t) \ .  
\end{equation}By lemma \ref{monodromy}, the holonomy of its orbits is   
\begin{equation}\label{holonomyfocusfocus}
[{hol_{\nabla^\omega}(\gamma)}](x_1,x_2,y_1,y_2)=\e^{i 2\pi (x_1y_2-x_2y_1)} \ .
\end{equation}\par

Now, given any $\boldsymbol\alpha\in S_{P|_W}^1(L|_W)$, using the unitary section $s$, it can be written as 
\begin{equation}
\boldsymbol\alpha=\alpha\otimes s=(\alpha_1\ud x_1 + \alpha_2\ud x_2 + \alpha_3\ud y_1 + \alpha_4\ud y_2)\otimes s 
\end{equation}and 
\begin{eqnarray}
\imath_{X_2}\phi_t^*\alpha|_{(x_1,x_2,y_1,y_2)}&=&[\alpha(\dot{\gamma}(t))]\circ\phi_t(x_1,x_2,y_1,y_2) \nonumber \\
&=&\alpha_1\circ\phi_t(x_1,x_2,y_1,y_2)(-x_1\sin t+x_2\cos t) \nonumber \\ 
&& +\alpha_2\circ\phi_t(x_1,x_2,y_1,y_2)(-x_2\sin t-x_1\cos t) \nonumber \\ 
&& +\alpha_3\circ\phi_t(x_1,x_2,y_1,y_2)(-y_1\sin t+y_2\cos t) \nonumber \\ 
&& +\alpha_4\circ\phi_t(x_1,x_2,y_1,y_2)(-y_2\sin t-y_1\cos t) \ .
\end{eqnarray}\par

Therefore, using 
\begin{equation}
A(t,p)=x_2\alpha_1\circ\phi_t(p)-x_1\alpha_2\circ\phi_t(p)+y_2\alpha_3\circ\phi_t(p)-y_1\alpha_4\circ\phi_t(p)
\end{equation}and
\begin{equation}
B(t,p)=x_1\alpha_1\circ\phi_t(p)+x_2\alpha_2\circ\phi_t(p)+y_1\alpha_3\circ\phi_t(p)+y_2\alpha_4\circ\phi_t(p) \ ,
\end{equation}the expression $\boldsymbol{J}_X(\boldsymbol\alpha)$ at a point $p=(x_1,x_2,y_1,y_2)$ near the singular point is:
\begin{eqnarray}
\left[\boldsymbol{J}_X(\boldsymbol\alpha)\right](p)&=&\left(\integral{0}{2\pi}{A(t,p)\e^{-it(x_1y_2-x_2y_1)}\cos t}{t}\right. 
\nonumber \\ &&\left.-\integral{0}{2\pi}{B(t,p)\e^{-it(x_1y_2-x_2y_1)}\sin t}{t}\right)s \ . 
\end{eqnarray}

The following upper bound proves that the expression above is zero over the set $\{(x_1,x_2,y_1,y_2)\in\CC^2; \ x_1y_2-x_2y_1=0\}$ (the points where ${hol_{\nabla^\omega}(\gamma)}=1$):
\begin{eqnarray}
|\boldsymbol{J}_X(\boldsymbol\alpha)|&\leq&\left|\max_{t\in[0,2\pi]}A(t,p)\right|\left|\integral{0}{2\pi}{\e^{-it(x_1y_2-x_2y_1)}\cos t}{t}\right| 
\nonumber \\ &&+\left|\max_{t\in[0,2\pi]}B(t,p)\right|\left|\integral{0}{2\pi}{\e^{-it(x_1y_2-x_2y_1)}\sin t}{t}\right| \nonumber \\
&=&\left|\max_{t\in[0,2\pi]}A(t,p)\right|\left|\fraction{(x_1y_2-x_2y_1)(\e^{-i2\pi(x_1y_2-x_2y_1)}-1)}{(x_1y_2-x_2y_1)^2-1}\right| 
\nonumber \\ &&+\left|\max_{t\in[0,2\pi]}B(t,p)\right|\left|\fraction{\e^{-i2\pi(x_1y_2-x_2y_1)}-1}{(x_1y_2-x_2y_1)^2-1}\right| \ .
\end{eqnarray}\par

The proof of lemma \ref{poincaeli2} works verbatim if $N$ is replaced by the local model $W$ describe here. Hence, this can be interpreted as a proof of the Poincaré lemma needed for the proof of theorem \ref{fineresolution1} when the real distribution has focus-focus singularities. 

\begin{teo}\label{poincafuck}
If $p\in M$ is a nondegenerate critical point of Williamson type \linebreak $(k_e,k_h,k_f)$, with $k_f\geq 1$ and $k_e+k_h+2k_f\leq n$ (it does not need to be a rank zero critical point, but it has to have a focus-focus component), for an integrable system $F:M\to\CR^n$ on a prequantisable symplectic manifold $(M,\omega)$, with polarisation induced by the moment map: the cohomology groups $H^k({S_P}^\bullet(L))$ vanish for $k\geq 0$ in a sufficiently small neighbourhood of $p$.
\end{teo}


\subsection{Neighbourhood of a Liouville fibre}

\hspace{1.5em}The Liouville theorem for integrable systems provides a symplectic normal form for a neighbourhood of a regular fibre. What follows is the computation of the quantisation of that model.\par
Let $M=\CR^n\times(\CR^{n-k}\times\mathds{T}^k)$ and $0\leq k\leq n$, where $\mathds{T}^k\cong\CR^k/2\pi\CZ^k$, with coordinates $(x_1,\dots,x_n,y_1,\dots,y_{n-k},\dots,y_n)$ and symplectic form $\omega=\somatorio{j=1}{n}\ud x_j\wedge\ud y_j$. It admits as a polarisation $P=\left\langle\depp{}{y_1},\dots,\depp{}{y_n}\right\rangle_{C^\infty(M)}$, and since $(M,\omega)$ is an exact symplectic manifold, it also admits as a prequantum line bundle $L=\CC\times M$ with connexion 1-form $\Theta=\somatorio{j=1}{n}x_j\ud y_j$, with respect to the unitary section $\exp\left(i\somatorio{j=1}{n}x_j\right)$.\par

The next lemma computes the contributions to geometric quantisation for each trivialising neighbourhood of a lagrangian fibre bundle.\par

\begin{lem}\label{prop555}  $H^{k+l}({S_P}^\bullet(L))=\{0\}$ for all $l\neq 0$ and 
\begin{equation}
H^k({S_P}^\bullet(L))\cong\left\{\begin{array}{l}\displaystyle\bigoplus_{m\in\CZ^k}C^\infty(\CR^{n-k}) \ , \ \text{if} \ k\neq n \\
\displaystyle\bigoplus_{m\in\CZ^k}\CC \ , \ \text{if} \ k=n\end{array}\right. \ .
\end{equation}
\end{lem}
\textit{Proof:} Supposing $k\neq n$, when $M$ is written as $(\CR\times S^1)^k\times\CC^{n-k}$ it becomes clear that the use of theorem \ref{kunnethregular} (more precisely, corollary \ref{kunnethregular2}) $k$ times reduces the problem of computing the quantisation of $M$ to the computation of the quantisation of $\CC^{n-k}$: which by proposition \ref{Qplano} is just $C^\infty(\CR^{n-k})$. If $k=n$, one just need to apply theorem \ref{kunnethregular} $n-1$ times, and, then, proposition \ref{Qcilindro} to conclude.\hfill $\blacksquare$\par\vspace{0.5em}


\subsection{Neighbourhood of an elliptic fibre}\label{ellipticfibre}

\hspace{1.5em}Toric, or locally toric, manifolds also have a normal form for a neighbourhood of its fibres, even if they are singular: Zung \cite{Zu1} attributes this normal form to Dufour and Molino, and Eliasson. The model and its quantisation are described below. 

For $0\leq k\leq n$, let $(M=\CR^{n+k}\times\mathds{T}^{n-k},\omega=\somatorio{j=1}{n}\ud x_j\wedge\ud y_j)$, with coordinates $(x_1,\dots,x_n,y_1,\dots,y_k,\dots,y_n)$, and $F:M\to\CR^n$ be a nondegenerate integrable system of elliptic type, i.e. 
\begin{equation}
F(x_1,\dots,x_n,y_1,\dots,y_n)=(x_1^2+y_1^2,\dots,x_k^2+y_k^2,x_{k+1},\dots,x_n) \ . 
\end{equation}

The polarisation in this case is 
\begin{equation}
P=\left\langle-y_1\depp{}{x_1}+x_1\depp{}{y_1},\dots,-y_k\depp{}{x_k}+x_k\depp{}{y_k},\depp{}{y_{k+1}},\dots,\depp{}{y_n}\right\rangle_{C^\infty(M)} \ ,
\end{equation}and since $(M,\omega)$ is an exact symplectic manifold, it also admits as a prequantum line bundle $L=\CC\times M$ with connexion 1-form $\Theta=\frac{1}{2}\somatorio{j=1}{k}(x_j\ud y_j-y_j\ud x_j)+\somatorio{j=k+1}{n}x_j\ud y_j$, with respect to the unitary section 
\begin{equation*}
\exp\left(i\somatorio{j=1}{k}(x_j^2+y_j^2)+i\somatorio{j=k+1}{n}x_j\right) \ .
\end{equation*}

\begin{prop}\label{prop556}  $\mathcal{Q}(M)\cong\CC^{b_s}$, where $b_s$ is the number of nonsingular Bohr-Sommerfeld fibres.
\end{prop}
\textit{Proof:} One can first use proposition \ref{kunnethelliptic} $k$ times, corollary \ref{kunnethregular2} $n-k-1$ times, and, finally, proposition \ref{Qcilindro}.\hfill $\blacksquare$\par\vspace{0.5em}

It is important to notice that, if the case of $b_s=0$ ($x_1^2+y_1^2,\dots,x_k^2+y_k^2 < 1$) was considered, the previous proof would give that all cohomology groups vanish when $k\neq 0$. This implies, as a corollary, that Bohr-Sommerfeld fibres of elliptic type give no contribution to geometric quantisation.


\subsection{Focus-focus contribution to geometric quantisation}\label{fuckfuck}

\hspace{1.5em}Let $F=(f_1,f_2):M\to\CR^2$ be an integrable system on a prequantasible $(M,\omega)$, with a nondegenerate focus-focus singular fibre $\ell_{ff}$ which is Bohr-Sommerfeld. In \cite{Zu2} it is demonstrated the existence of a neighbourhood of $\ell_{ff}$ over which $f_2$ is a moment map for a hamiltonian $S^1$-action. And proposition \ref{rawaction} asserts that there are no nonzero flat sections near a nondegenerate focus-focus singular fibre.

\begin{coro}
In the neighbourhood of $\ell_{ff}$ over which a hamiltonian $S^1$-action is defined, there exists a neighbourhood $V$ containing only $\ell_{ff}$ as a Bohr-Sommerfeld fibre such that $\check{H}^0(V;\mathcal{J}\big{|}_V)=\{0\}$.
\end{coro} 

Concerning how focus-focus fibres behave under geometric quantisation, this partially answers a conjecture raised by Hamilton and Miranda: stating that there is no contribution coming from focus-focus singularities to geometric quantisation.

Believing that the conjecture is true, one could try to use proposition \ref{rawpotential} to prove it for the neighbourhood $V$. The first obstacle is that $\{{hol_{\nabla^\omega}(\gamma)}=1\}$ is not a submanifold, and one needs to prove that $\boldsymbol{J}_X$ is the null operator over the points where $\{{hol_{\nabla^\omega}(\gamma)}=1\}$ (which, looking at the proof of theorem \ref{poincafuck}, seems to be the case). Another approach would be to prove only the exactness of $\boldsymbol{J}_X$ and check out convergence over the singular points of $\{{hol_{\nabla^\omega}(\gamma)}=1\}$. 


\section{Applications II: global computations}\label{applications}


\subsection{Lagrangian fibre bundles}

\hspace{1.5em}In \cite{Sni} \'{S}niatycki studies the case when the polarisation is a lagrangian fibration. He uses a resolution for the sheaf and proves the vanishing of the groups $\check{H}^l(M;\mathcal{J})$, for $l\neq k$: $k$ being the rank of the fundamental group of a fibre. 

\begin{teo}[\'{S}niatycki]\label{sniteo} If the base space $N$ is a manifold and the natural projection $\mathcal{F}:M\to N$ is a lagrangian fibration, then $\mathcal{Q}(M)=\check{H}^k(M;\mathcal{J})$, and $\check{H}^k(M;\mathcal{J})\cong\check{H}^0(\ell_{BS};\mathcal{J}|_{\ell_{BS}})$, where $\ell_{BS}\subset M$ is the union of all Bohr-Sommerfeld fibres.
\end{teo}

A slightly different proof of his theorem is given here when $k\neq 0$. When $k=0$ there is no symplectic circle action and the techniques presented in chapter \ref{circleraw} are of no use; wherefore, apart from the presentation, the proof is the same as the original one and is omitted.\par

Any atlas of the base space satisfies that the projection $\mathcal{F}:M\to N$ on each open set $V$ of the atlas is a moment map. Assuming $\mathrm{dim}(M)=2n$, if $\chi:V\to\CR^n$ is a coordinate system over $V$, $F:=\chi\circ\mathcal{F}\big{|}_{\mathcal{F}^{-1}(V)}:\mathcal{F}^{-1}(V)\to\CR^n$ is an integrable system because each $f_j:=\mathrm{pr}_j\circ F$ is constant along the fibres of $\mathcal{F}$, $\ud f_j=0$ along them. In other words, $\ud f_j$ annihilates vector fields tangent to the fibres and the hamiltonian vector fields of the others $f_j$'s are, indeed, tangent to the fibres: $\{f_i,f_j\}_\omega=0$.  

The open sets $\mathcal{F}^{-1}(V)$ are just the model in lemma \ref{prop555} with a fixed number of Bohr-Sommerfeld fibres; thus, the quantisation of it is just a sum of copies of $\CC$, or $C^{\infty}(\CR^{n-k})$, depending on the value of $k$, for each Bohr-Sommerfeld fibre.\par 

Assuming that $k\neq 0$, so that theorem \ref{integerlattice2} can be used, the atlas can ---and it will--- be chosen in such a way that no Bohr-Sommerfeld fibre is contained in more than one of the open sets $\mathcal{F}^{-1}(V)$. In particular, if $V$ and $W$ are two open sets of the atlas such that $V\cap W\neq\emptyset$, then $\mathcal{F}^{-1}(V)\cap\mathcal{F}^{-1}(W)$ has no Bohr-Sommerfeld fibre. Proposition \ref{bsleaf} implies that one of the periodic hamiltonian vector fields of the components of $F$ has orbits with nontrivial holonomy over $\mathcal{F}^{-1}(V)\cap\mathcal{F}^{-1}(W)$; thus, by proposition \ref{rawpotential}, its quantisation is just the trivial vector space $\{0\}$. This means that a Mayer-Vietoris argument works for the cover $\{\mathcal{F}^{-1}(V)\}$ of $M$, and this finishes the proof for $k\neq 0$.\par

\begin{rema} \'{S}niatycki works with the prequantum line bundle twisted by a bundle of half forms normal to the polarisation, and here the result is presented for the nontwisted prequantum line bundle. As it was mentioned before, the techniques used here apply to any complex line bundle admitting a flat connexion along the polarisation: the only difference being that the Bohr-Sommerfeld fibres may not be the same.  
\end{rema}

\begin{exam}
The Kodaira-Thurston manifold is an example of a lagrangian bundle (there is a description of it in \cite{Ha}). Moreover, it is a symplectic manifold which is not Kähler: it gives at least one reason for developing a trully symplectic geometric quantisation apparatus.\hfill $\Diamond$\par\vspace{0.5em}
\end{exam}


\subsection{Locally toric manifolds}

\hspace{1.5em}Hamilton \cite{Ha} has shown, via \v{C}ech approach, that \'{S}niatycki's theorem holds for locally toric manifolds and that the elliptic fibres give no contribution to the quantisation.

\begin{teo}[Hamilton]\label{hateo} For $M$ a $2n$-dimensional compact symplectic manifold equipped with a locally toric singular lagrangian fibration: 
\begin{equation}
\mathcal{Q}(M)=\check{H}^n(M;\mathcal{J})\cong\displaystyle\bigoplus_{p\in BS_r}\CC \ ,
\end{equation}$BS_r$ being the set of the regular Bohr-Sommerfeld fibres.
\end{teo}

\begin{rema} Regarding metaplectic correction, contrary to \'{S}niatycki's, Hamilton's result does not include a twisted prequantum line bundle. Using the framework described in this article, it is straightforward to twist the prequantum line bundle by a bundle of half forms normal to the polarisation and achieve the same result ---only noticing that the Bohr-Sommerfeld fibres may not be the same.  
\end{rema}

The previous reasoning used for the fibre bundle case works in this singular setting. This provides a proof for Hamilton's theorem via a de Rham approach.\par  

A locally toric singular lagrangian fibration on a symplectic manifold $(M,\omega)$ is a surjective map $\mathcal{F}:M\to N$, where $N$ is a topological space such that for every point in $N$ there exist an open neighbourhood $V$ and a homeomorphism $\chi:V\to U\subset\{z\in\CR^k \ ; \ z\geq 0\}\times\CR^{n-k}$ satisfying that $(\mathcal{F}^{-1}(V),\omega\big{|}_{\mathcal{F}^{-1}(V)},\chi\circ\mathcal{F}\big{|}_{\mathcal{F}^{-1}(V)})$ is an integrable system symplectomorphic to an open subset of the model of proposition \ref{prop556}.\par

Hence, by definition, the open sets $\mathcal{F}^{-1}(V)$ are just the model in proposition \ref{prop556} with a fixed number of Bohr-Sommerfeld fibres; thus, the quantisation of it is just a sum of copies of $\CC$, or $\{0\}$, depending on the fibre dimension, for each Bohr-Sommerfeld fibre.\par 

Choosing an open cover for $N$ in such a way that no Bohr-Sommerfeld fibre is contained in more than one of the open sets $\mathcal{F}^{-1}(V)$ (theorem \ref{integerlattice2} allows one to make this choice), if $V$ and $W$ are two open sets of the atlas such that $V\cap W\neq\emptyset$, then $\mathcal{F}^{-1}(V)\cap\mathcal{F}^{-1}(W)$ has no Bohr-Sommerfeld fibre. Proposition \ref{bsleaf} implies that one of the periodic hamiltonian vector fields of the components of the integrable system has orbits with nontrivial holonomy over $\mathcal{F}^{-1}(V)\cap\mathcal{F}^{-1}(W)$; wherefore, by proposition \ref{rawpotential}, its quantisation is just the trivial vector space $\{0\}$. This means that a Mayer-Vietoris argument works for the cover $\{\mathcal{F}^{-1}(V)\}$ of $M$, and this finishes the proof.\par  


\subsection{Almost toric manifolds}

\hspace{1.5em}As it was seen from the quantisation of lagrangian fibrations and locally toric manifolds, quantisation of neighbourhoods of Bohr-Sommerfeld fibres computes the quantisation of the whole manifold. Consequently, if one knows how to compute the higher cohomology groups for a neighbourhood of a focus-focus fibre, one is able to compute the quantisation for the almost toric case using the factorisation tools (corollary \ref{kunnethregular2} and proposition \ref{kunnethelliptic}) and proceeding like the lagrangian bundle and locally toric cases.

For example, if $M$ is a $4$-dimensional compact almost toric manifold, $BS_r$ and $BS_{ff}$ are the image of the regular, respectively focus-focus, Bohr-Sommerfeld fibres by the moment map, and $V$ neighbourhoods of focus-focus fibres admitting a hamiltonian $S^1$-action:
\begin{equation}
\mathcal{Q}(M)\cong\left(\bigoplus_{p\in BS_r}\CC\right)\oplus\left(\bigoplus_{p\in BS_{ff}}\check{H}^1(V;\mathcal{J}\big{|}_V)\oplus\check{H}^2(V;\mathcal{J}\big{|}_V)\right) \ . 
\end{equation}\par
    


\end{document}